\documentstyle[12pt,dina4,amsfont12]{article} 
 
\input psfig.sty 
 
\begin{document} 
 
 
\newcommand{\pl}[2]{\frac{\partial#1}{\partial#2}} 
\newcommand{\Ta}{{\cal A}} 
\newcommand{\Tb}{{\cal B}} 
\newcommand{\Te}{{\cal E}} 
\newcommand{\bu}{{\bf u} } 
\newcommand{\p}{\partial} 
\newcommand{\og}{\omega} 
\newcommand{\Og}{\Omega} 
\newcommand{\fl}[2]{\frac{#1}{#2}} 
\newcommand{\dt}{\delta} 
\newcommand{\tm}{\times} 
\newcommand{\sm}{\setminus} 
\newcommand{\nn}{\nonumber} 
\newcommand{\ap}{\alpha} 
\newcommand{\bt}{\beta} 
\newcommand{\ld}{\lambda} 
\newcommand{\Gm}{\Gamma} 
\newcommand{\gm}{\gamma} 
\newcommand{\vp}{\varphi} 
\newcommand{\tht}{\theta} 
\newcommand{\ift}{\infty} 
\newcommand{\vep}{\varepsilon} 
\newcommand{\ep}{\epsilon} 
\newcommand{\kp}{\kappa} 
\newcommand{\Dt}{\Delta} 
\newcommand{\Sg}{\Sigma} 
\newcommand{\fa}{\forall} 
\newcommand{\sg}{\sigma} 
\newcommand{\ept}{\emptyset} 
\newcommand{\btd}{\nabla} 
\newcommand{\btu}{\Delta} 
\newcommand{\tg}{\triangle} 
\newcommand{\Th}{{\cal T}_h} 
\newcommand{\ged}{\qquad \Box} 
\newcommand{\bgv}{\Bigg\vert} 
\renewcommand{\theequation}{\arabic{section}.\arabic{equation}} 
\newcommand{\be}{\begin{equation}} 
\newcommand{\ee}{\end{equation}} 
\newcommand{\ba}{\begin{array}} 
\newcommand{\ea}{\end{array}} 
\newcommand{\bea}{\begin{eqnarray}} 
\newcommand{\eea}{\end{eqnarray}} 
\newcommand{\beas}{\begin{eqnarray*}} 
\newcommand{\eeas}{\end{eqnarray*}} 
\newcommand{\dpm}{\displaystyle} 
\newtheorem{theorem}{Theorem}[section] 
\newtheorem{lemma}{Lemma}[section] 
\newtheorem{remark}{Remark}[section] 
\newcommand{\Gmu}{\Gm_{_U}} 
\newcommand{\Gml}{\Gm_{_L}} 
\newcommand{\Gme}{\Gm_e} 
\newcommand{\Gmi}{\Gm_i} 
\newcommand{\lN}{{_N}} 
\newcommand{\tld}[1]{\~{#1}} 
\newcommand{\td}[1]{\tilde{#1}} 
\newcommand{\um}{\mu} 
\newcommand{\bp}{{\bf \Psi} } 
\newcommand{\bx}{{\bf x} } 
 
\title{An explicit unconditionally stable numerical method 
 for solving damped 
nonlinear Schr\"{o}dinger equations with a focusing nonlinearity} 
\author{ \ 
{\it Weizhu Bao} 
\thanks{Email address:  bao@cz3.nus.edu.sg.}\\ 
Department of Computational Science\\ 
National University of Singapore, Singapore 117543\\ 
\\ 
{\it Dieter Jaksch } 
\thanks{Email address: dieter.jaksch@physics.oxford.ac.uk.}\\ 
Institut f\"{u}r Theoretische Physik, Universit\"at Innsbruck,\\ 
A--6020 Innsbruck, Austria.\\ 
} 
 
\date{} 
\maketitle 
 
\begin{abstract} 
This paper introduces an extension of the time-splitting 
sine-spectral (TSSP) method for solving damped focusing nonlinear 
Schr\"{o}dinger equations (NLS). The method is explicit, 
unconditionally stable and time transversal invariant. Moreover, 
it preserves the exact decay rate for the normalization of the 
wave function if linear damping terms are added to the NLS. 
Extensive numerical tests are presented for cubic focusing 
nonlinear Schr\"{o}dinger equations in 2d with a linear, cubic or 
a quintic damping term. Our numerical results show that quintic 
or cubic damping always arrests blowup, while linear damping can 
arrest blowup only when the damping parameter $\dt$ is larger than 
a threshold value $\dt_{\rm th}$. We note that our method can also 
be applied to solve the 3d Gross-Pitaevskii equation with a 
quintic damping term to model the dynamics of a collapsing and 
exploding Bose-Einstein condensate (BEC). 
\end{abstract} 
 
\bigskip 
 
 
\bigskip 
 
 
\section{Introduction}\label{si} 
\setcounter{equation}{0} 
 
Since the first experimental realization of Bose-Einstein 
condensation (BEC) in dilute weakly interacting gases the 
nonlinear Schr\"{o}dinger equation (NLS) has been used extensively 
to describe the single particle properties of BECs. The results 
obtained by solving the NLS showed excellent agreement with most 
of the experiments (for a review see 
\cite{Anglin,Dalfovo,Cornell}). In fact, up to now there have been 
very few experiments in ultracold dilute bosonic gases which could 
not be described properly by using theoretical methods based on 
the NLS \cite{NGPEExp,NGPETheo}. 
 
Recent experiments by Donley {\em et al.} \cite{Donley} provide 
new experimental results for checking the validity of describing a 
BEC by using the NLS in the case of attractive interactions 
(focusing nonlinearity) in 3d. Since the particle density might 
become very large in the case of attractive interactions inelastic 
collisions become important and cannot be neglected. These 
inelastic collisions are assumed to be accounted for by adding 
damping terms to the NLS. Two particle inelastic processes are 
taken into account by a cubic damping term while three particle 
inelastic collisions are described by a quintic damping term. 
Collisions with the background gas and feeding of the condensate 
can be studied by adding linear damping terms. One of the major 
theoretical challenges in comparing results obtained in the 
experiment with theoretical results is to find reliable methods 
for solving the NLS with a focusing nonlinearity and damping terms 
in the parameter regime where the experiments are performed. 
 
The aim of this paper is to extend the time-splitting 
sine-spectral method (TSSP) for solving the focusing NLS with 
additional damping terms and to present extensive numerical tests. 
The comparison of our numerical results with the experimental 
results obtained for a collapsing BEC \cite{Donley} will be 
presented elsewhere \cite{BaoBEC}. 
 
We consider the NLS \cite{Bao3,Sulem} 
\begin{eqnarray} \label{sdge} 
&&i\; \psi_t=-\fl{1}{2}\;\btu \psi+ V(\bx)\; \psi -\bt 
|\psi|^{2\sg}\psi, \qquad t>0, \qquad \bx\in {\Bbb R}^d, \\ 
\label{sdgi} &&\psi(\bx,t=0)=\psi_0(\bx), \qquad \bx\in {\Bbb 
R}^d, 
\end{eqnarray} 
with $\sg>0$ a positive constant, where $\sg=1$ corresponds to a cubic 
nonlinearity and $\sg=2$ corresponds to a quintic nonlinearity, 
$V({\bf x})$ is a real-valued potential whose shape is determined 
by the type of system under investigation, and $\bt$ 
positive/negative corresponds to the focusing/defocusing NLS. In 
BEC, where (\ref{sdge}) is also known as the Gross-Pitaevskii 
equation (GPE) \cite{Pit}, $\psi$ is the macroscopic wave function 
of the condensate, $t$ is time, ${\bf x}$ is the spatial 
coordinate and $V({\bf x})$ is a trapping potential which usually 
is harmonic and can thus be written as $V({\bf x})=\fl{1}{2} 
\left(\gm_1^2 x_1^2 +\cdots+\gm_d^2 x_d^2\right)$ with $\gm_1, 
\cdots, \gm_d\ge 0$. Two important invariants of (\ref{sdge}) are 
the {\bf normalization of the wave function} 
\begin{equation} \label{mass} 
N(t)=\int_{{\Bbb R}^d}\; |\psi({\bf x}, t)|^2\; d{\bf x}, 
\qquad t\ge 0 
\end{equation} 
and the {\bf energy} 
\begin{equation} \label{energy} 
E(t) = \int_{{\Bbb R}^d}\; \left[\fl{1}{2}|\btd \psi({\bf x}, 
t)|^2 
+V(\bx)|\psi(\bx,t)|^2-\fl{\bt}{\sg+1}|\psi(\bx,t)|^{2\sg+2}\right]\; 
d\bx, \quad t\ge 0. 
\end{equation} 
 
From the theory for the local existence of the solution of 
(\ref{sdge}), it is well known that if $\|\psi(\cdot,t)\|_{H^1}$ 
is bounded, the solution exists for all $t$ \cite{Sulem}. As a 
result, when the NLS is defocusing ($\bt<0$), conservation of 
energy implies that $\int_{{\Bbb R}^d}\;|\btd \psi({\bf x}, t)|^2 
\;d\bx$ is bounded and the solution exists globally. On the other 
hand, if the NLS is focusing ($\bt>0$) at critical ($\sg d=2$) or 
supercritical ($\sg d>2$) dimensions and for an initial energy 
$E(0)<0$, the solutions of (\ref{sdge}) can self-focus and become 
singular in finite time, i.e. there exists a time $t_*<\ift$ such 
that \cite{Sulem} 
\[ 
\lim_{t\to t_*} |\btd \psi|_{L^2}=\ift \qquad 
{\rm and} \qquad  \lim_{t\to t_*} |\psi|_{L^\ift}=\ift. 
\] 
However, the physical quantities modeled by $\psi$ do not become 
infinite which implies that the validity of (\ref{sdge}) breaks 
down near the singularity. Additional physical mechanisms, which 
were initially small, become important near the singular point and 
prevent the formation of the singularity. In BEC the particle 
density $|\psi|^2$ becomes large close to the critical point and 
inelastic collisions between particles which are negligible for 
small densities become important. Therefore a small damping 
(absorption) term is introduced into the NLS (\ref{sdge}) which 
describes inelastic processes. We are interested in the cases 
where these damping mechanisms are important and, therefore, 
restrict ourselves to the case of focusing nonlinearities $\bt>0$, 
where $\bt$ may also be time dependent. We consider the following 
damped nonlinear Schr\"{o}dinger equation: 
\begin{eqnarray} \label{sdged} 
&&i\; \psi_t=-\fl{1}{2}\;\btu \psi+ V(\bx)\; \psi -\bt 
|\psi|^{2\sg}\psi-i\; g(|\psi|^2)\psi, 
\qquad t>0, \quad \bx\in {\Bbb R}^d, \\ 
\label{sdgid} 
&&\psi(\bx,t=0)=\psi_0(\bx), \qquad \bx\in {\Bbb 
R}^d, 
\end{eqnarray} 
where $g(\rho)\ge 0$ for $\rho=|\psi|^2\ge 0$ is a real-valued 
monotonically increasing function. 
 
The general form of (\ref{sdged}) covers many damped NLS arising 
in various different applications. In BEC, for example, when 
$g(\rho)\equiv 0$, (\ref{sdged}) reduces to the usual GPE 
(\ref{sdge}); a linear damping term $g(\rho)\equiv \dt$ with 
$\dt>0$ describes inelastic collisions with the background gas; 
cubic damping $g(\rho)=\dt_1 \bt \rho$ with $\dt_1>0$ corresponds 
to two-body loss \cite{Saito,Roberts}; and a quintic damping term 
of the form $g(\rho)=\dt_2 \bt^2 \rho^2$ with $\dt_2>0$ adds 
three-body loss to the GPE (\ref{sdge}) \cite{Saito,Roberts}. It's 
easy to see that the decay of the normalization according to 
(\ref{sdged}) due to damping is given by 
\be \label{normNt} 
N^\prime(t)=\fl{d}{dt} \int_{{\Bbb R}^d}\; |\psi({\bf x}, t)|^2\; 
d{\bf x} = -2\int_{{\Bbb R}^d}\; g(|\psi({\bf x}, t)|^2)|\psi({\bf 
x}, t)|^2\; d{\bf x} \le 0, \quad t >0. 
\ee 
Particularly, if 
$g(\rho)\equiv \dt$ with $\dt>0$, the normalization is given by 
\be \label{drnt} 
N(t)=\int_{{\Bbb R}^d}\; |\psi({\bf x}, 
t)|^2\; d{\bf x}= e^{-2\dt\; t} N(0)=e^{-2\dt\; t}\; \int_{{\Bbb 
R}^d}\; |\psi_0({\bf x})|^2\; d{\bf x}, \quad t\ge 0. 
\ee 
 
There has been a series of recent studies which deals with the 
analysis and numerical solution of the damped NLS. Fibich 
\cite{Fibich} analyzed the effect of linear damping (absorption) 
on the critical self-focusing NLS, Tsutsumi \cite{Tsutsumi, 
TsutsumiM} studied the global solutions of the NLS with linear 
damping, the regularity of attractors and approximate inertial 
manifolds for a weakly damped NLS were given in Goubet 
\cite{Goubet,GoubetOl} and by Jolly {\rm et al.}~\cite{Jolly}. For 
numerically solving the linearly damped NLS Peranish 
\cite{Peranich} proposed a finite difference scheme and this 
method was revisited recently by Ciegis {\em et al.}~\cite{Ciegis} 
and Zhang {\em et al.}~\cite{Zhang}. Moebs \cite{Moebs} presented 
a multilevel method for weakly damped NLS and applied it to solve 
a stochastic weakly damped NLS in \cite{MoebsG}. Variable mesh 
difference schemes for the NLS with a linear damping term were 
used by Iyengar {\em et al.}~\cite{Iyengar}. 
 
Also the TSSP, which we will use in this paper to solve the damped 
NLS, was already successfully used for solving the Schr\"{o}dinger 
equation in the semiclassical regime and for describing 
Bose-Einstein condensation using the Gross-Pitaeskii equation by 
Bao {\em et al.} \cite{Bao1,BJP,Bao3}. The TSSP is explicit, 
unconditionally stable and time transversal invariant. Moreover, 
it gives the exact decay rate of the normalization when linear 
damping is applied to the NLS (i.e. $g(\rho)\equiv \dt$ with 
$\dt>0$ in (\ref{sdged})) and yields spectral accuracy for spatial 
derivatives and second-order accuracy for the time derivative. 
Thus this method is a very good candidate for solving the damped 
NLS, especially in 2d or 3d. We test the novel numerical method 
extensively in 2d. 
 
Finally, we want to emphasize that the NLS is also used in 
nonlinear optics, e.g., to describe the propagation of an intense 
laser beam through a medium with a Kerr nonlinearity 
\cite{FibichP, Sulem}. In nonlinear optics $\psi=\psi({\bf x},t)$ 
describes the electrical field amplitude, $t$ is the spatial 
coordinate in the direction of propagation, ${\bf x} =(x_1,\cdots, 
x_d)^T$ is the transverse spatial coordinate and $V({\bf x})$ is 
determined by the index of refraction. Nonlinear damping terms of 
the form $g(\rho)=\dt \bt^q \rho^q$ with $\dt, q>0$ correspond to 
multiphoton absorption processes \cite{Fibich}. 
 
The paper is organized as follows. In section \ref{smethod} we 
present the time-splitting sine-spectral approximation for the 
damped nonlinear Schr\"{o}dinger equation. In section \ref{sne} 
numerical tests are presented for the cubic focusing nonlinear 
Schr\"{o}dinger equation in 2d with a linear, cubic or quintic 
damping term. In section \ref{sc} some conclusions are drawn. 
 
\section{Time-splitting sine-spectral method} 
\label{smethod} 
\setcounter{equation}{0} 
 
In this section we present a time-splitting sine-spectral (TSSP) 
method for solving the problem (\ref{sdged}), (\ref{sdgid}) with 
homogeneous periodic boundary conditions. For simplicity of 
notation we shall introduce the method for the case of one spatial 
dimension $(d=1)$. Generalizations to $d>1$ are straightforward 
for tensor product grids and the results remain valid without 
modifications. For $d=1$, the problem becomes 
\begin{eqnarray} \label{sdge1d} 
&&i\;\psi_t=- \fl{1}{2} \psi_{xx}+ V(x)\psi 
-\bt|\psi|^{2\sg}\psi -i\, g(|\psi|^2)\psi, 
\qquad a<x<b,\quad t>0, \qquad \\ 
\label{sdgi1d} &&\psi(x,t=0)=\psi_0(x), \quad a\le x\le b, \qquad 
\psi(a,t)=\psi(b,t)={\bf 0}, \quad t\ge 0. 
\end{eqnarray} 
 
\subsection{General damping term} 
 
We choose the spatial mesh size $h=\btu x>0$ with $h=(b-a)/M$ and 
$M$ an even positive integer, the time step is given by $k=\btu 
t>0$ and define grid points and time steps by 
\[ 
x_j:=a+j\;h, \qquad t_n := n\; k, \qquad j=0,1,\cdots, M, \qquad 
n=0,1,2,\cdots 
\] 
Let $\psi^{n}_j$ be the numerical approximation of $\psi(x_j,t_n)$ 
and $\psi^{n}$ the solution vector at time $t=t_n=nk$ with 
components $\psi_j^{n}$. 
 
 From time $t=t_n$ to time $t=t_{n+1}$, the damped nonlinear 
Schr\"{o}dinger equation (\ref{sdge1d}) is solved in two steps. 
One solves 
\begin{equation} \label{fstep} 
i\; \psi_t=- \fl{1}{2} \psi_{xx}, 
\end{equation} 
for one time step, followed by solving 
\begin{equation} 
\label{sstep} i\; \psi_t(x,t)= V(x)\psi(x,t)- \bt 
|\psi(x,t)|^{2\sg} \psi(x,t)-i\, g(|\psi(x,t)|^2)\psi(x,t), 
\end{equation} 
again for the same time step. Equation (\ref{fstep}) is 
discretized in space by the sine-spectral method and integrated in 
time {\it exactly}. For $t\in[t_n,t_{n+1}]$, multiplying the ODE 
(\ref{sstep}) by $\overline{\psi(x,t)}$, the conjugate of 
$\psi(x,t)$, one obtains 
\begin{equation} \label{sstepa} 
i\;\psi_t(x,t)\overline{\psi(x,t)}= V(x)|\psi(x,t)|^2- \bt 
|\psi(x,t)|^{2\sg+2} -i\, g(|\psi(x,t)|^2)|\psi(x,t)|^2. 
\end{equation} 
Subtracting the conjugate of Eq.~(\ref{sstepa}) from 
Eq.~(\ref{sstepa}) and multiplying by $-i$ one obtains 
\begin{equation} 
\label{sstepc} 
\fl{d}{dt}|\psi(x,t)|^2=\overline{\psi_t(x,t)}\psi(x,t)+ 
\psi_t(x,t)\overline{\psi(x,t)} =-2 g(|\psi(x,t)|^2)|\psi(x,t)|^2. 
\end{equation} 
Let 
\begin{equation} \label{frho} 
f(s)=\int \fl{1}{s\; g(s)}\;ds, \qquad h(s,\tau)=\left\{\ba{ll} 
f^{-1}\left(f(s)-2\tau\right), &s>0, \ \tau\ge0, \\ 
0, &s=0, \ \tau\ge0. \ea\right. \\ 
\end{equation} 
Then, if $g(s)\ge0$ for $s\ge 0$, we find 
\begin{equation} 
\label{positive} 0\le h(s,\tau)\le s, \qquad \hbox{for}\quad s\ge 
0,\quad  \tau\ge 0 
\end{equation} 
and the solution of the ODE (\ref{sstepc}) can be expressed as 
(with $\tau=t-t_n$) 
\begin{eqnarray} \label{rhot} 0\le\rho(t)&=&\rho(t_n+\tau) 
:=|\psi(x,t)|^2=h\left(|\psi(x,t_n)|^2, t-t_n\right):= 
h\left(\rho(t_n),\tau\right) \nn\\ 
&\le&\rho(t_n)=|\psi(x,t_n)|^2, \qquad t_n\le t\le t_{n+1}. 
\end{eqnarray} 
Combining Eq.~(\ref{rhot}) and Eq.~(\ref{sstep}) we obtain 
\begin{eqnarray} \label{sstepd} 
i\;\psi_t(x,t)&=&V(x)\psi(x,t)- 
\bt \left[h\left(|\psi(x,t_n)|^2, t-t_n\right)\right]^{\sg} \psi(x,t)\nn\\ 
&&-i\, g\left(h\left(|\psi(x,t_n)|^2,t-t_n\right)\right)\psi(x,t), 
\qquad t_n\le t\le t_{n+1}. 
\end{eqnarray} 
Integrating (\ref{sstepd}) from $t_n$ to $t$, we find 
\begin{eqnarray} 
\label{spsi} \psi(x,t)&=& 
\exp\left\{i\left[-V(x)(t-t_n)+G\left(|\psi(x,t_n)|^2,t-t_n\right)\right] 
-F\left(|\psi(x,t_n)|^2,t-t_n\right)\right\} \nn\\ 
&&\quad \tm\; \psi(x,t_n), \qquad t_n\le t\le t_{n+1}, 
\end{eqnarray} 
where we have defined 
\begin{equation} 
\label{FG} 
F(s,r)=\int_0^r g(h(s,\tau))\;d\tau\ge0, 
\quad 
G(s,r)=\int_0^r \bt\; \left[h(s,\tau)\right]^{\sg}\;d\tau. 
\end{equation} 
To find the time evolution between $t=t_n$ and $t=t_{n+1}$, we 
combine the splitting steps via the standard second-order Strang 
splitting (TSSP) for solving the damped nonlinear Schr\"{o}dinger 
equation (\ref{sdge1d}). In detail, the steps for obtaining 
$\psi^{n+1}_j$ from $\psi^{n}_j$ are given by 
\begin{eqnarray} 
&&\psi^*_j=\exp\left\{-F\left(|\psi_j^n|^2,k/2\right)+i 
\left[-V(x_j)k/2+G\left(|\psi_j^n|^2,k/2\right)\right]\right\} 
\;\psi_j^{n}, \nn\\ 
\label{schmg} 
&&\psi_j^{**}=\sum_{l=1}^{M-1} e^{-i 
k\mu_l^2/2}\;\widehat{\psi}^*_l \; \sin(\mu_l(x_j-a)),\qquad 
j=1,2,\cdots,M-1,\\ 
&&\psi^{n+1}_j 
=\exp\left\{-F\left(|\psi_j^{**}|^2,k/2\right)+i 
\left[-V(x_j)k/2+G\left(|\psi_j^{**}|^2,k/2\right)\right]\right\} 
\;\psi_j^{**}, \nn 
\end{eqnarray} 
where $ \widehat{U}_l$ are the sine-transform coefficients of a 
complex vector $U=(U_0,U_1, \cdots, U_M)$ with $U_0=U_M={\bf 0}$ 
which are defined as 
\begin{equation} 
\label{Fourc1} \mu_l=\fl{\pi l}{b-a}, \quad \widehat{U}_l= 
\fl{2}{M}\sum_{j=1}^{M-1} U_j\;\sin(\mu_l (x_j-a)), \ 
l=1,2,\cdots, M-1, 
\end{equation} 
where 
\begin{equation} 
\label{init1} 
\psi^{0}_j=\psi(x_j,0)=\psi_0(x_j), \qquad j=0,1,2,\cdots,M. 
\end{equation} 
Note that the only time discretization error of TSSP is the 
splitting error, which is second order in $k$ if the integrals in 
(\ref{frho}) and (\ref{FG}) can be evaluated analytically. 
 
\subsection{Most frequently used damping terms} 
 
In this subsection we present explicit formulae for using TSSP when 
solving the NLS with those damping terms most frequently appearing 
in BEC and nonlinear optics. 
 
\bigskip 
 
\noindent {\bf Case I} NLS with a linear damping term. We choose 
$g(\rho)\equiv \dt$ with $\dt>0$ in (\ref{sdged}). In BEC this 
damping terms describes inelastic collisions of condensate 
particles with the background gas. From (\ref{frho}), we find 
\begin{equation} 
\label{frhol} 
f(s)=\int \fl{1}{\dt s}ds=\fl{1}{\dt}\ln s \qquad 
{\rm and}\qquad h(s,\tau)=e^{-2\dt \tau }\;s. 
\end{equation} 
Substituting (\ref{frhol}) into (\ref{rhot}) and (\ref{FG}), we obtain 
\begin{eqnarray} \label{FG11l} 
 &&\rho(t)= e^{-2\dt (t-t_n) }\; |\psi(x,t_n)|^2, \quad t_n\le t\le t_{n+1},\\ 
\label{FG12l} 
 &&F(s,r)=\dt r, \\ 
\label{FG13l} 
 &&G(s,r)=\fl{\bt s^\sg}{2\dt\sg}\left(1-e^{-2\dt \sg r}\right). 
\end{eqnarray} 
Substituting (\ref{FG12l}) and (\ref{FG13l}) into (\ref{schmg}), 
we get the following second-order time-splitting sine-spectral 
steps for the NLS with a linear damping term 
\begin{eqnarray} 
&&\psi^*_j=\exp\left\{-k\dt/2+i 
\left[-V(x_j)k/2+\bt |\psi_j^n|^{2\sg}\left(1-e^{-\dt\sg k} 
\right)/(2\dt\sg)\right]\right\} 
\;\psi_j^{n},   \nn\\ 
\label{schmld} 
&&\psi_j^{**}=\sum_{l=1}^{M-1} 
  e^{-i k\mu_l^2/2}\;\widehat{\psi}^*_l\;\sin(\mu_l(x_j-a)), 
    \qquad j=1,2,\cdots,M-1,\\ 
&&\psi^{n+1}_j=\exp\left\{-k\dt/2+i 
\left[-V(x_j)k/2+\bt |\psi_j^{**}|^{2\sg}\left(1-e^{-\dt\sg k} 
\right)/(2\dt\sg)\right]\right\} 
\;\psi_j^{**}.  \nn 
\end{eqnarray} 
 
\noindent {\bf Case II} NLS with a damping term of the form 
$g(\rho)= \dt\bt^q \rho^q $, where $\dt,\; q>0$ in (\ref{sdged}). 
For $q=1$ ($q=2$) we obtain the damping term describing two 
(three) particle inelastic collisions in BEC. From (\ref{frho}) we 
get 
\begin{equation} 
\label{frhogd} 
f(s)=\int \fl{1}{\dt \bt^q s^{q+1}}ds=-\fl{1}{q\dt\bt^q s^q} \qquad 
{\rm and}\qquad h(s,\tau)=\fl{s}{\left(1+2q\dt\tau \bt^q s^q \right)^{1/q}}. 
\end{equation} 
Substituting (\ref{frhogd}) into (\ref{rhot}) and (\ref{FG}), we obtain 
\begin{eqnarray} 
\label{FG41l} 
 &&\rho(t)= \fl{|\psi(x,t_n)|^2}{\left[1+2q\dt\bt^q(t-t_n) 
 |\psi(x,t_n)|^{2q} \right]^{1/q}}, \quad t_n\le t\le t_{n+1}, \qquad \\ 
\label{FG42l} 
 &&F(s,r)= \fl{1}{2q}\ln\left(1+2q\dt r \bt^q s^q\right), \qquad\quad  \\ 
\label{FG43l} 
 &&G(s,r)=\dpm\left\{\ba{ll} \dpm\fl{\bt^{1-q}}{2\dt q}\ln\left(1+2q\dt 
 r \bt^q s^q \right) &q=\sg, \\ 
 \dpm\fl{\bt^{1-q}s^{\sg-q}\left[-1+(1+2q\dt r 
 \bt^q s^q)^{(q-\sg)/q}\right]} 
 {2\dt (q-\sg)}  &\sg\ne q.\quad \\ 
 \ea\right. 
\end{eqnarray} 
Substituting (\ref{FG42l}) and (\ref{FG43l}) into (\ref{schmg}), 
we get the following second-order time-splitting sine-spectral 
method for the NLS 
{\footnotesize 
\begin{eqnarray} 
\label{schmgd} 
&&\psi^*_j=\dpm\left\{\ba{ll} 
\dpm\fl{\exp\left\{i 
\left[-V(x_j)k/2+\bt^{1-q}\ln\left(1+\dt q k \bt^q 
|\psi_j^n|^{2q}\right)/(2\dt q) 
\right]\right\}}{\left(1+q\dt k\bt^q |\psi_j^n|^{2q}\right)^{1/2q}} 
\;\psi_j^{n},  &\sg=q,\\ 
\ &\ \\ 
\dpm\fl{\exp\left\{i 
\left[-\fl{V(x_j)k}{2}+\fl{\bt^{1-q}|\psi_j^n|^{2\sg-2q}}{2\dt(q-\sg)} 
\left(-1+\left(1+\dt q k \bt^q 
|\psi_j^n|^{2q}\right)^{\fl{q-\sg}{q}}\right) 
\right]\right\}}{\left(1+q\dt k\bt^q |\psi_j^n|^{2q}\right)^{1/2q}} 
\;\psi_j^{n},  &\sg\ne q,\\ 
\ea\right. \nn\\ 
&&\psi_j^{**}=\sum_{l=1}^{M-1} 
  e^{-i k\mu_l^2/2}\;\widehat{\psi}^*_l\;\sin(\mu_l(x_j-a)), 
    \qquad j=1,2,\cdots,M-1,\\ 
&&\psi^{n+1}_j=\dpm\left\{\ba{ll} 
\dpm\fl{\exp\left\{i 
\left[-V(x_j)k/2+\bt^{1-q}\ln\left(1+\dt q k \bt^q 
|\psi_j^{**}|^{2q}\right)/(2\dt q) 
\right]\right\}}{\left(1+q\dt k\bt^q |\psi_j^{**}|^{2q}\right)^{1/2q}} 
\;\psi_j^{**},  &\sg=q,\\ 
\ &\ \\ 
\dpm\fl{\exp\left\{i 
\left[-\fl{V(x_j)k}{2}+\fl{\bt^{1-q}|\psi_j^{**}|^{2\sg-2q}}{2\dt(q-\sg)} 
\left(-1+\left(1+\dt q k \bt^q 
|\psi_j^{**}|^{2q}\right)^{\fl{q-\sg}{q}}\right) 
\right]\right\} }{\left(1+q\dt k\bt^q |\psi_j^{**}|^{2q}\right)^{1/2q}} 
\;\psi_j^{**},  &\sg\ne q.\\ 
\ea\right. \nn 
\end{eqnarray} 
} 
 
\bigskip 
 
\noindent {\bf Case III} Focusing cubic NLS with a damping term 
that accounts for two-body and three-body loss in a BEC 
\cite{Saito}, i.e., we choose $\sg=1$, $g(\rho)= \dt_1\bt 
\rho+\dt_2\bt^2\rho^2 $ with $\dt_1,\; \dt_2>0$, in (\ref{sdged}). 
Using (\ref{frho}), we get 
\begin{eqnarray} 
\label{frhogdd} f(s)&=&\left\{\ba{ll} -\fl{1}{\dt_1\bt 
s}+\fl{\dt_2}{\dt_1^2}\ln\left(\dt_2\bt+\dt_1/s\right), &s>0,\\ 
 0, &s=0. \ea\right. 
\end{eqnarray} 
Substituting (\ref{frho}) into (\ref{FG}) and changing the 
variable of integration we obtain 
\begin{eqnarray} 
\label{FGd42l} 
 F(s,r)&=& \int_0^r g\left( f^{-1}(f(s)-2\tau)\right)\;d\tau 
\stackrel{\tau=(f(s)-f(h))/2}{=}\int_s^{h(s,r)} 
-\fl{1}{2} g(h) f^\prime(h)\;dh \nn\\ 
&=&\int_s^{h(s,r)}-\fl{1}{2h}\;dh=\left\{\ba{ll} 
-\fl{1}{2}\ln\left(h(s,r)/s\right), &s>0,\\ 
0 &s=0;\\ 
\ea\right. 
\end{eqnarray} 
where $h(s,r)$ is the solution of 
\begin{equation} 
\label{shsr} 
f(s)-f(h(s,r))=2r, \qquad \hbox{for any}\ r>0, 
\end{equation} 
with $f$ given in (\ref{frhogdd}). Similarly we find 
\begin{eqnarray} 
\label{FGd43l} 
 G(s,r)&=&\int_s^{h(s,r)}-\fl{\bt }{2g(h)}\;dh 
=\left\{\ba{ll} -\fl{1}{2\dt_1} \ln\fl{h(s,r)(\dt_1+\dt_2\bt s)} 
{s(\dt_1+\dt_2\bt h(s,r))}, &s>0, \\ 
0, &s=0. \\ 
\ea\right. 
\end{eqnarray} 
Substituting (\ref{FGd42l}) and (\ref{FGd43l}) into (\ref{schmg}) 
we get the following second-order time-splitting sine-spectral 
steps for the NLS with a combination of cubic and quintic damping 
terms {\scriptsize 
\begin{eqnarray} 
\label{schmgdd} 
&&\psi^*_j=\dpm\left\{\ba{ll} 
\dpm\fl{\sqrt{h(|\psi_j^n|^2,k/2)}}{|\psi_j^n|} 
\exp\left\{i 
\left[-\fl{V(x_j)k}{2}-\fl{1}{2\dt_1}\ln 
\fl{h(|\psi_j^n|^2,k/2)(\dt_1+\dt_2\bt|\psi_j^n|^2)} 
{|\psi_j^n|^2(\dt_1+\dt_2\bt h(|\psi_j^n|^2,k/2))}\right]\right\} 
\;\psi_j^{n},  &\psi_j^n\ne 0,\\ 
\ &\ \\ 
0,  &\psi_j^n=0,\\ 
\ea\right. \nn\\ 
&&\psi_j^{**}=\sum_{l=1}^{M-1} 
  e^{-i k\mu_l^2/2}\;\widehat{\psi}^*_l\;\sin(\mu_l(x_j-a)), 
    \qquad j=1,2,\cdots,M-1,\\ 
&&\psi^{n+1}_j=\dpm\left\{\ba{ll} 
\dpm\fl{\sqrt{h(|\psi_j^{**}|^2,k/2)}}{|\psi_j^{**}|} 
\exp\left\{i 
\left[-\fl{V(x_j)k}{2}-\fl{1}{2\dt_1}\ln 
\fl{h(|\psi_j^{**}|^2,k/2)(\dt_1+\dt_2\bt|\psi_j^{**}|^2)} 
{|\psi_j^{**}|^2(\dt_1+\dt_2\bt h(|\psi_j^{**}|^2,k/2))}\right]\right\} 
\;\psi_j^{**},  &\psi_j^{**}\ne 0,\\ 
\ &\ \\ 
0,  &\psi_j^{**}=0.\\ 
\ea\right. \nn 
\end{eqnarray} 
} 
 
\begin{remark} 
As demonstrated in this subsection, the integrals in (\ref{frho}) 
and (\ref{FG}) can be evaluated {\bf analytically} for the damping 
terms which most frequently appear in physical applications. If 
the integrals in (\ref{frho}) or (\ref{FG}) can not be evaluated 
analytically or the inverse of $f$ in (\ref{frho}) can not be 
expressed explicitly, e.g., if $g(\rho)$ in (\ref{sdged}) is not a 
polynomial, one can solve the following ODE numerically by either 
second- or fourth-order Runge-Kutta method \beas \label{ODEe1} 
&&\fl{d h(t)}{dt}=-2g(h(t))\; h(t), \qquad 0\le t\le k/2, \\ 
\label{ODEe2} &&h(0)=s, \eeas to get $h(s,k/2)$ for any given 
$s>0$ and set $h(s,k/2)=0$ for $s=0$. By changing the variable of 
integration in (\ref{FG}), see detail in (\ref{FGd42l}) and 
(\ref{FGd43l}), the first integral in (\ref{FG}), i.e. $F(s,k/2)$, 
can be evaluated exactly (see detail in (\ref{FGd42l})), and the 
second integral in (\ref{FG}), i.e. $G(s,k/2)=\dpm 
\int_s^{h(s,k/2)} -\fl{\bt h^{\sg-1}}{2g(h)}\;dh$, can be 
evaluated numerically by using a numerical quadrature, e.g., the 
trapezoidal rule or Simpson's rule. 
\end{remark} 
 
The scheme TSSP is explicit and is unconditionally stable as we 
will demonstrate in the next subsection. Another main advantage of 
the time-splitting method is its time transversal invariance, which 
also holds for the NLS and the damped NLS themselves. If a 
constant $\ap$ is added to the potential $V$, then the discrete 
wave functions $\psi_j^{\vep,n+1}$ obtained from TSSP get 
multiplied by the phase factor $e^{-i\ap(n+1)k}$, which leaves the 
discrete normalization unchanged. This property does not hold for 
finite difference schemes.

\begin{remark} 
For the focusing cubic NLS with a quintic damping term describing 
three-body recombination loss and an additional feeding term for 
the BEC \cite{Kagan} we choose $\sg=1$, $g(\rho)= -\dt_1+ 
\dt_2\bt^2 \rho^2 $ with $\dt_1,\dt_2>0$ in (\ref{sdged}). The 
idea of constructing the TSSP is also applicable to this case 
although we could not prove that it is unconditonally stable due 
to the feeding term. Inserting the above feeding term into 
(\ref{frho}), we get 
\begin{eqnarray} 
\label{frhoqdd} f(s)&=&\left\{\ba{ll} 
\fl{1}{2\dt_1}\ln \left|\dt_2\bt^2 -\dt_1/s^2\right|, &s>0,\\ 
0, &s=0.\\ 
\ea\right. 
\end{eqnarray} 
Inserting (\ref{frhoqdd}) into (\ref{rhot}), we find 
\begin{equation} 
\label{hstgd} h(s,\tau)=\fl{s \sqrt{\dt_1}} {\sqrt{ \dt_1 
e^{-4\tau \dt_1}+(1-e^{-4\tau\dt_1})\dt_2\bt^2 s^2}}, 
\end{equation} 
and substituting (\ref{hstgd}) into (\ref{rhot}) and (\ref{FG}), 
we obtain 
\begin{eqnarray} 
\label{FGg31l} 
 &&\rho(t)=\fl{|\psi(x,t_n)|^2\sqrt{\dt_1}} {\sqrt{ 
\dt_1 e^{-4\tau \dt_1}+(1-e^{-4\tau\dt_1})\dt_2\bt^2 
|\psi(x,t_n)|^4}}, 
\ t_n\le t\le t_{n+1}, \qquad \quad \\ 
\label{FGg32l} 
 &&F(s,r)=-\dt_1 r+\fl{1}{4}\ln\left[1+ \dt_2\bt^2 
s^2(e^{4\dt_1 r}-1)/\dt_1\right], 
\qquad \\ 
\label{FGg33l} 
 &&G(s,r)=\fl{1}{2\sqrt{\dt_1\dt_2}}\ln \fl{\bt 
s\sqrt{\dt_2}e^{2r \dt_1}+\sqrt{\dt_1+ \dt_2 \bt^2 
s^2\left(e^{4r\dt_1}-1\right)}}{\sqrt{\dt_1}+\bt s\sqrt{\dt_2}}. 
\end{eqnarray} 
Inserting (\ref{FGg32l}) and (\ref{FGg33l}) into (\ref{schmg}), we 
get the following second-order time-splitting sine-spectral steps 
for the NLS with a quintic damping term and a feeding term 
{\footnotesize 
\begin{eqnarray} 
\label{schdgdd} 
&&\psi^*_j=\dpm\fl{e^{k\dt_1/2}\exp\left[i\left(-\fl{V(x_j)k}{2}+ 
\fl{1}{2\sqrt{\dt_1\dt_2}}\ln 
\fl{\bt |\psi_j^n|^2\sqrt{\dt_2}e^{k \dt_1}+\sqrt{\dt_1+ 
\dt_2 \bt^2 |\psi_j^n|^4\left(e^{2k\dt_1}-1\right)}} 
{\sqrt{\dt_1}+\bt |\psi_j^n|^2\sqrt{\dt_2}}\right)\right]} 
{\left[1+\dt_2\bt^2|\psi_j^n|^4(e^{2k\dt_1}-1)/\dt_1\right]^{1/4}} 
\;\psi_j^{n},\nn\\ 
&&\psi_j^{**}=\sum_{l=1}^{M-1} 
  e^{-i k\mu_l^2/2}\;\widehat{\psi}^*_l\;\sin(\mu_l(x_j-a)), 
    \qquad j=1,2,\cdots,M-1,\\ 
&&\psi^{n+1}_j=\dpm\fl{e^{k\dt_1/2}\exp\left[i\left(-\fl{V(x_j)k}{2}+ 
\fl{1}{2\sqrt{\dt_1\dt_2}}\ln 
\fl{\bt |\psi_j^{**}|^2\sqrt{\dt_2}e^{k \dt_1}+\sqrt{\dt_1+ 
\dt_2 \bt^2 |\psi_j^{**}|^4\left(e^{2k\dt_1}-1\right)}} 
{\sqrt{\dt_1}+\bt |\psi_j^{**}|^2\sqrt{\dt_2}}\right)\right]} 
{\left[1+\dt_2\bt^2|\psi_j^{**}|^4(e^{2k\dt_1}-1)/\dt_1\right]^{1/4}} 
\;\psi_j^{**}.  \nn 
\qquad \quad 
\end{eqnarray} } 
\end{remark}

\begin{remark} 
The scheme TSSP (\ref{schmg}) can easily be extended for solving the 
complex Ginzburg-Landau equation (CGL) \cite{FibchLevy,Mielke} \be 
\label{cgl} i\; \psi_t=-\left(1-i\;\vep\right)\btu \psi - 
|\psi|^{2}\psi -i\left(\dt_2|\psi|^2-\dt_1\right)\psi, \ee where 
$\vep$, $\dt_1$ and $\dt_2$ are positive constants. The idea of 
construncting the TSSP for the damped NLS is also applicaple to 
the CGL provided that we solve \be \label{spcgl} i\; 
\psi_t=-\left(1-i\;\vep\right)\btu \psi, \ee in the first step 
instead of (\ref{fstep}). Inserting $\sg=1$, $\bt=1$  and 
$g(\rho)=\dt_2 \rho-\dt_1$ with $\dt_1$, $\dt_2>0$ into 
(\ref{sdged}) and using (\ref{frho}) we get 
\begin{eqnarray} 
\label{frhocgl} 
f(s)&=&\left\{\ba{ll} 
\fl{1}{\dt_1}\ln \left|\dt_2 -\dt_1/s\right|, &s>0,\\ 
0, &s=0.\\ 
\ea\right. 
\end{eqnarray} 
Inserting (\ref{frhocgl}) into (\ref{frho}) we find 
\begin{equation} 
\label{hstgdcgl} 
h(s,\tau)=\fl{s \dt_1} {s\dt_2\left(1-e^{-2\tau \dt_1}\right) 
+\dt_1e^{-2\tau \dt_1} }, 
\end{equation} 
and substituting (\ref{hstgdcgl}) into (\ref{rhot}) and (\ref{FG}) 
we obtain 
\begin{eqnarray} 
\label{FGg31lcgl} 
 &&\rho(t)=\fl{\dt_1\;|\psi(x,t_n)|^2} {\dt_2\;|\psi(x,t_n)|^2 
\left(1-e^{-2\tau \dt_1}\right)+\dt_1e^{-2\tau \dt_1} }, 
\ t_n\le t\le t_{n+1}, \qquad \quad \\ 
\label{FGg32lcgl} 
 &&F(s,r)=-\fl{1}{2}\ln \fl{\dt_1}{s\dt_2+\left(\dt_1-s\dt_2\right) 
e^{-2r \dt_1}}, \qquad \\ 
\label{FGg33lcgl} 
 &&G(s,r)=\fl{1}{2\dt_2}\ln \fl{\dt_1-s\dt_2+s\dt_2 e^{2r\dt_1}}{\dt_1}. 
\end{eqnarray} 
Inserting (\ref{FGg32lcgl}) and (\ref{FGg33lcgl}) into (\ref{schmg}), we 
get the following second-order time-splitting sine-spectral steps 
for the CGL (\ref{cgl}) 
{\footnotesize 
\begin{eqnarray} 
\label{schdgcgl} 
&&\psi^*_j=\sqrt{\fl{\dt_1}{\dt_2|\psi_j^n|^2+\left(\dt_1-\dt_2 
|\psi_j^n|^2\right) 
e^{-k\dt_1}} }\exp\left[\fl{i}{2\dt_2}\ln 
\fl{\dt_1-\dt_2|\psi_j^n|^2+\dt_2|\psi_j^n|^2e^{k\dt_1}}{\dt_1}\right] 
\;\psi_j^{n},  \nn\\ 
&&\psi_j^{**}=\sum_{l=1}^{M-1} 
  e^{-(\vep+i)k\mu_l^2}\;\widehat{\psi}^*_l\;\sin(\mu_l(x_j-a)), 
    \qquad j=1,2,\cdots,M-1,\\ 
&&\psi^{n+1}_j=\sqrt{\fl{\dt_1}{\dt_2|\psi_j^{**}|^2+\left(\dt_1-\dt_2 
|\psi_j^{**}|^2\right) 
e^{-k\dt_1}} }\exp\left[\fl{i}{2\dt_2}\ln 
\fl{\dt_1-\dt_2|\psi_j^{**}|^2+\dt_2|\psi_j^{**}|^2e^{k\dt_1}}{\dt_1}\right] 
\;\psi_j^{**},  \nn 
\end{eqnarray} } 
\end{remark} 
 
\begin{remark} 
If the homogeneous peridic boundary conditions in (\ref{sdgi1d}) 
are replaced by the periodic boundary conditions \be \label{gpc} 
\psi(a,t)=\psi(b,t), \qquad \psi_x(a,t)=\psi_x(b,t), \qquad t\ge0, 
\ee the TSSP scheme (\ref{schmg}) still works provided that one 
replaces the sine-series in (\ref{schmg}) by a Fourier-series 
\cite{Bao,BJP,Bao3}. 
\end{remark}

\subsection{Stability and decay rate} 
 
Let $U=(U_0, U_1, \cdots, U_{M})^T$ with $U_0=U_M={\bf 0}$ and 
$\|\cdot\|_{l^2}$ be the usual discrete $l^2$-norm on the interval 
$(a,b)$, i.e., 
\begin{equation} 
\label{norm} 
\|U\|_{l^2}=\sqrt{\fl{b-a}{M}\sum_{j=1}^{M-1} |U_j|^2}. 
\end{equation} 
For the {\it stability} of the time-splitting sine-spectral 
approximations  TSSP (\ref{schmg}), we have the following lemma, 
which shows that the total normalization does not increase. 
 
\begin{lemma}\label{stability} 
The time-splitting sine-spectral schemes (TSSP) (\ref{schmg}) are 
unconditionally stable if $g(s)\ge0$ for $s\ge 0$. In fact, for 
every mesh size $h>0$ and time step $k>0$, 
\begin{equation} 
\label{stabu} 
\|\psi^{n+1}\|_{l^2}\le \|\psi^{n}\|_{l^2}\le 
\|\psi^{0}\|_{l^2}=\|\psi_0\|_{l^2}, 
 \qquad n=0,1,2,\cdots 
\end{equation} 
Furthermore, when a linear damping term is used in (\ref{sdged}), 
i.e., we choose $g(\rho)\equiv \dt$ with $\dt>0$, the decay rate 
of the normalization satisfies 
\begin{equation} 
\label{decayr} 
\|\psi^{n}\|_{l^2} = e^{-2\dt t_n}\|\psi^{0}\|_{l^2} 
=e^{-2\dt t_n}\|\psi_0\|_{l^2}, \qquad n=1,2,\cdots 
\end{equation} 
In fact, (\ref{decayr}) is a discretized version of the 
decay rate of the normalization $N(t)$ in (\ref{drnt}). 
\end{lemma} 
 
\noindent {Proof:}\ We combine (\ref{schmg}), (\ref{Fourc1}), 
(\ref{norm}) and note that $F(s,\tau)\ge 0$ for $s\ge0$ and 
$\tau\ge0$, to obtain 
\begin{eqnarray} 
\label{stabSP2} 
\lefteqn{\fl{1}{b-a}\|\psi^{n+1}\|_{l^2}^2=\fl{1}{M}\sum_{j=1}^{M-1} 
  \left|\psi^{n+1}_j\right|^2}\nn\\[2mm] 
&=&\fl{1}{M}\sum_{j=1}^{M-1} 
\exp\left[-2F\left(|\psi_j^{**}|^2,k/2\right)\right] 
\;\left|\psi_j^{**}\right|^2 
  \le \fl{1}{M}\sum_{j=1}^{M-1}\left|\psi_j^{**}\right|^2 \nn\\ 
  &=&\fl{1}{M}\sum_{j=1}^{M-1}\left| 
\sum_{l=1}^{M-1} 
  e^{-i k\mu_l^2/2}\;\widehat{\psi}^*_l\;\sin(\mu_l(x_j-a))\right|^2 
  =\fl{1}{2}\sum_{l=1}^{M-1}\left| 
   e^{-i k\mu_l^2/2}\; \hat{\psi}^{*}_l\right|^2 
=\fl{1}{2}\sum_{l=1}^{M-1}\left| 
    \hat{\psi}_l^{*}\right|^2 \nn\\ 
&=&\fl{1}{2}\sum_{l=1}^{M-1}\left|\fl{2}{M} 
   \sum_{j=1}^{M-1} \psi_j^{*}\; \sin(\mu_l (x_j-a))\right|^2 
=\fl{1}{M} 
   \sum_{j=1}^{M-1} \left|\psi_j^{*}\right|^2\nn\\ 
&=&\fl{1}{M} 
   \sum_{j=1}^{M-1} \exp\left[-2F\left(|\psi_j^n|^2,k/2\right)\right] 
\left|\psi_j^{n}\right|^2 
\le \fl{1}{M} 
   \sum_{j=1}^{M-1} \left|\psi_j^{n}\right|^2 \nn\\ 
&=&\fl{1}{b-a}\|\psi^{n}\|^2_{l^2}. 
\end{eqnarray} 
Here, we used the identity 
\begin{equation} 
\label{iden1} 
\sum_{j=1}^{M-1} \sin\left(\fl{\pi r\; j}{M}\right) 
\; \sin\left(\fl{\pi s\; j}{M}\right)=\left\{\ba{ll} 
0, &\ r-s\ne 2mM, \\ 
M/2, &\ r-s=2mM, r\ne 2nM, 
\ea\right. \quad m,n\ {\rm integer}. 
\end{equation} 
 
When a linear damping term is added to the NLS (\ref{sdged}), the 
equality (\ref{decayr}) follows from the above proof, 
Eq.~(\ref{FG12l}), and 
\[   \sum_{j=1}^{M-1} \exp\left[-2F\left(|\psi_j^n|^2,k/2\right)\right] 
\left|\psi_j^{n}\right|^2 = 
   \sum_{j=1}^{M-1} e^{-\dt k} \left|\psi_j^{n}\right|^2 
=e^{-\dt k}   \sum_{j=1}^{M-1}  \left|\psi_j^{n}\right|^2. 
\]

\section{Numerical examples}\label{sne} 
\setcounter{equation}{0} 
 
In this section we present numerical tests of the TSSP 
(\ref{schmg}) for solving a focusing cubic NLS appearing in 
nonlinear optics \cite{FibichP, Sulem} and for the Gross-Pitaeskii 
equation in BEC \cite{Bao3} in 2d with a linear, a cubic, or a 
quintic damping term. In our computations, the initial condition 
(\ref{sdgi}) is always chosen such that $|\psi_0(\bx)|$ decays to 
zero sufficiently fast as $|{\bf x}|\to\ift$. We choose an 
appropriately large rectangle $[a,b]\tm[c,d]$ in 2d to avoid that 
the homogeneous periodic boundary condition (\ref{sdgi1d}) 
introduce a significant (aliasing) error relative to the whole 
space problem. To quantify the numerical results of the GPE for a 
BEC, we define the condensate widths along the $x$, $y$ and 
$z$-axis by 
\[ 
\sigma_\ap^2=\langle \ap^2\rangle=\fl{1}{N(t)}\;\int_{{\Bbb R}^d}\; \ap^2 
|\psi({\bf x},t)|^2\; d{\bf x}, \qquad \hbox{with}\quad 
\ap=x,\ y,\ {\rm or}\ z. 
\] 
 
\bigskip 
 
{\bf Example 1} Solution of the 2d  damped  focusing cubic 
nonlinear Schr\"{o}dinger equation. We choose $d=2$, $\sg=1$ and 
$V(x,y)\equiv0$ in (\ref{sdged}) and present computations for 
three different damping terms ($\dt>0$): 
 
\bigskip 
 
\noindent {\it I. A linear damping term}, i.e. we choose $g(\rho)\equiv 
\dt$. 
 
\smallskip 
\noindent {\it II. A cubic damping term}, i.e. we choose $g(\rho)\equiv 
\dt \bt \rho$. 
 
\smallskip 
\noindent {\it III. A quintic damping term}, i.e. we choose $g(\rho)\equiv 
\dt \bt^2 \rho^2$. 
 
\bigskip 
 
The initial condition (\ref{sdgid}) is taken to be 
\be 
\label{initdata} 
\psi(x,y,0)=\psi_0(x,y)=\fl{\gm_y^{1/4}}{\sqrt{\pi \vep}}\; 
e^{-(x^2+\gm_y y^2)/2\vep}, \qquad (x,y)\in {\Bbb R}^2. 
\ee 
We assume $\gm_y=2$, $\vep=0.2$ and $\bt=8$ in (\ref{sdged}) such 
that $E(0)=-0.751582<0$ in (\ref{energy}). We solve the NLS on the 
square $[-16,16]^2$, i.e., $a=c=-16$ and $b=d=16$ with mesh size 
$h=\fl{1}{32}$, time step $k=0.0002$ and homogeneous periodic 
boundary conditions along the boundary of the square. We compare 
the effect of changing the damping parameter $\dt$ in the three 
different cases I, II and III. 
 
Figure 1 shows the surface plot of the density $|\psi(x,y,t)|^2$ 
at time $t=1.25$ with $\dt=0.5$; plots of the normalization, energy and 
central density $|\psi(0,0,t)|^2$ as functions of time 
with $\dt=0.5,\ 0.3$ and $\dt=0$ (no damping) for case I. 
Figure 2  shows similar results for case II and Figure 3 
for case III. Furthermore Figure 4 shows contour plots of the 
density $|\psi|^2$ at  different times for case III with $\dt=0.01$. 
 
  In the numerically computations, a blowup is detected either from 
the plot of the central density $|\psi(0,0,t)|^2$ which at the blowup 
shows a very sharp spike with a peak value that increases when the 
mesh size $h$ decreases, or from the plot of the energy $E(t)$ which 
has a very sharp spike with negative values at the blowup. 
In fact, the method TSSP (\ref{schmg}) aims to capture the solution 
of damped NLS without blowup, i.e. physical revelant solution. 
If one wants to capture the blowup rate of NLS, we refer to 
\cite{LPW,WP}. 
 
\bigskip 
 
 From the numerical results we find the following conditions for 
arresting a blowup of the wave function with initial energy 
$E(0)<0$. (1) For linear damping the blowup is arrested if the 
damping parameter is bigger than a certain threshold value which 
we find to be $\dt_{\rm th} \approx 0.461$ by numerical 
experiments. As shown in Fig.~1b blowup is arrested for $\dt = 0.5 
> \dt_{\rm th}$ while the wave function blows up for $\dt < 
\dt_{\rm th}$ as can be seen from Fig.~1c\&d where we have chosen 
$\dt= 0.3 < \dt_{\rm th}$ and $\dt= 0 < \dt_{\rm th}$, 
respectively. The time at which the blowup of the wave function 
happens, however, increases with increasing $\dt$ (cf.~Fig.~1c\&d). 
(2) For a cubic damping term with 
$\dt>0$ the blowup of the wave function is always arrested 
 (cf.~Fig.~2). (3) The above observation (2) also holds for 
 a quintic damping term (cf.~Fig.~3). 
 
  For linear damping, we also test the dependence of the threshold value 
of the damping parameter $\dt_{\rm th}$ on $\bt$ and the initial data. 
First we take $\gm_y=2$ and  $\vep=0.2$ in (\ref{initdata}). Table 1 shows 
the threshold values $\dt_{\rm th}$ for different $\bt$ in (\ref{sdged}), 
and $E(0)$ represents the initial energy. Then we choose $\bt=16$ in 
(\ref{sdged}) and $\gm_y=2$ in (\ref{initdata}). Table 2 displays the 
threshold values $\dt_{\rm th}$ for different values of $\vep$ 
in (\ref{initdata}). 
 
\begin{table}[htbp] 
\begin{center} 
\begin{tabular}{cccccc}\hline 
    &$\bt=8$ &$\bt=16$ &$\bt=32$ 
   &$\bt=64$ &$\bt=128$\\ 
\hline 
$E(0)$  &$-0.7516$ &$-5.253$ &$-14.256$ &$-32.263$ &$-68.275$\\ 
$\dt_{\rm th}$ &$0.461$ &$3.655$ &$10.35$ &$22.15$ &$40.05$\\ 
\hline 
\end{tabular} 
\end{center} 
Table 1: Dependence of $\dt_{\rm th}$ on $\bt$ for 
$\gm_y=2$ and $\vep=0.2$ in (\ref{initdata}). 
\end{table} 
 
\begin{table}[htbp] 
\begin{center} 
\begin{tabular}{cccccc}\hline 
    &$\vep=0.8$ &$\vep=0.4$ &$\vep=0.2$ 
  &$\vep=0.1$  &$\vep=0.05$ \\ 
\hline 
$E(0)$  &$-1.3133$ &$-2.6266$ &$-5.2532$ &$-10.506$ &$-21.013$\\ 
$\dt_{\rm th}$ &$0.895$ &$1.845$ &$3.655$ &$7.25$ &$14.55$\\ 
\hline 
\end{tabular} 
\end{center} 
Table 2: Dependence of $\dt_{\rm th}$ on $\vep$ in (\ref{initdata}) for 
$\bt=16$ in (\ref{sdged}) and $\gm_y=2$ in (\ref{initdata}). 
\end{table} 
 
 From Table 1 we find by a least square fitting, 
\[\dt_{\rm th}=-0.6930 E(0) \qquad \hbox{or}\qquad 
 \dt_{\rm th}=0.3872 \bt-2.4627.\] 
Similarly, from Table 2 we obtain 
\[\dt_{\rm th}=-0.6922 E(0).\] 
Based on this observation, we conclude that the threshold value 
of the linear damping parameter $\dt_{\rm th}$ depends linearly 
on the initial energy $E(0)$. 
 
\bigskip 
 
{\bf Example 2} Solution of the 2d  damped GPE with focusing 
nonlinearity. We choose $d=2$, $\sg=1$ and $V(x,y)= 
\fl{1}{2}(\gm_x^2 x^2+\gm_y^2 y^2)$ to be a harmonic oscillator 
potential with $\gm_x,\gm_y>0$ in (\ref{sdged}). Again, we present 
computations for the same three different damping terms in (\ref{sdged}) 
as those we studied in Example 1. 
 
We take $\gm_x=1$ and $\gm_y=4$. The initial condition 
(\ref{sdgid}) is assumed to be the ground-state solution of 
(\ref{sdged}) with $g(\rho)\equiv0$ (i.e.~undamped case) and 
$\bt=-40$. The cubic nonlinearity is ramped linearly from 
$\bt=-40$ (defocusing) to $\bt=50$ (focusing) during the time 
interval $[0,0.1]$ and afterwards kept constant. The absorption 
parameter was set to $\dt=0$ during the time interval $[0,0.1]$ 
and increased to a positive value $\dt>0$ afterwards. 
 
We solve the GPE on the rectangle $[-24,24]\tm[-6,6]$, i.e., for 
$a=-24$, $b=24$, $c=-6$ and $d=6$ with mesh size $h_x=\fl{3}{64}$, 
$h_y=\fl{3}{128}$, time step $k=0.0005$ and homogeneous periodic 
boundary conditions along the boundary of the rectangle. Again, we 
compare the effect of changing the damping parameter $\dt$ in the 
three different cases I, II and III. 
 
Figure 5 shows a surface plot of the density $|\psi(x,y,t)|^2$ at 
times $t=0$ (ground-state solution) and $t=2.8$ with $\dt=1.25$; 
normalization, energy and central density $|\psi(0,0,t)|^2$ as 
functions of time with $\dt=1.25$, $1.1$ and $0$ (no damping) for 
case I. Figure 6 shows similar results for case II and Figure 7 
for case III. Furthermore Figure 8 shows contour plots of the 
density $|\psi|^2$ at  different times for case III with $\dt=0.15$. 
 
\bigskip 
 
 From our numerical results we find that the observations (1)-(3) 
made for example 1 are still valid with the additional trapping 
potential. However, the value of $\dt_{\rm th}$, depends on $\bt$ 
(or initial energy $E(0)$) and we find 
$\dt_{\rm th}\approx 1.185$ for linear damping 
(cf.~Fig.~5). 
 
\subsection{Discussion} 
 
In this subsection we discuss our numerical results in terms of 
physical properties of a BEC described by the GPE. We concentrate 
on those cases where a collapse of the wave function is arrested 
since this collapse leads to unphysical processes like the 
negative peaks in the energy $E(t)$ shown in Figs.~1c\&d,5e\&f. 
 
The general form of the time evolution in example 1 is similar for 
all three cases. Initially the cloud of atoms contracts due to the 
attractive interaction between the particles. This contraction is 
accompanied by an increase in the energy due to particle loss 
which is most efficient in regions of high particle density. These 
regions are characterized by a negative local energy density 
leading to an increase in energy for each particle lost there. 
After the central particle density has reached a maximum the cloud 
starts to expand due to the kinetic energy gained by the particles 
during the contraction. Particles are emitted from the cloud in 
burst like pulses which can be seen in Figs.~4,8. Such bursts 
have also been seen in BEC experiments \cite{Donley}. The main 
differences between the three cases are the behavior of the energy 
and the number of particles as a function of time. In case I where 
we assumed a linear damping term the loss rate of particles from 
the condensate is independent of the shape of the condensate wave 
function. The energy decrease during the condensate expansion is 
determined by the loss of particles (cf.~Fig.~1b). In the cases of 
cubic and quintic damping the loss term only has a significant 
effect on the time evolution of the condensate during the 
contraction. When the condensate expands the density of particles 
is so low that the loss terms have only a very small effect and 
the energy $E(t)$ and the number of particles $N(t)$ remain almost 
constant (see Figs.~2c,3c\&d). 
 
In example 2 we add an additional trap potential which confines 
the BEC and assume a realistic scenario (described above) to prepare the 
condensate in the trap (cf. experiments by Donley {\em et 
al.} \cite{Donley}). We find that the initial process of turning 
on the attractive interactions between the particles leads to 
oscillations in the widths of the condensate \cite{Bao3} as can be 
seen from Figs.~5,6,7. However, neither the additional trap 
potential nor these oscillations significantly alter the behavior 
of the system compared to example 1 when the condensate is 
strongly contracted. Before and after this contraction some 
differences can be seen. By looking at Figs.~5,6 we find that 
the first minimum in $\sigma_y$ due to the oscillations of the 
condensate causes and increase in the central density and in the 
energy. For cubic and quintic damping this is accompanied by an 
increased particle loss. However, an arrested collapse of the wave 
function only happens when both $\sigma_x$ and $\sigma_y$ attain a 
minimum value due to the attractive interactions (cf.~Fig.~5d 
and Fig.~6b). We also note that the frequency of the oscillations 
after an arrested collapse has happened is not significantly 
influenced by the damping terms. The amplitude of these 
oscillations is, however, strongly dependent on $\dt$ and 
decreases with increasing $\dt$. Finally, we want to mention that 
a series of contractions and expansion of the condensate is 
possible. In Fig.~7b we find three contractions of the condensate 
where only the first one reaches a sufficiently high particle 
density to lead to an increase in energy while the next two 
contractions show a rather smooth decrease in energy and particle 
number. For a smaller quintic damping term we obtain two 
contractions of the condensate which increase the energy (see 
Fig.~7c). 
 
\section{Conclusions}\label{sc} 
\setcounter{equation}{0} 
 
We extended the explicit unconditionally stable second-order 
time-splitting sine-spectral (TSSP) method for solving damped 
focusing nonlinear Schr\"{o}dinger equations. We showed that this 
method is time transversal invariant and preserves the exact decay 
rate of the normalization for a linear damping of the NLS. 
Extensive numerical tests were presented for the cubic focusing 
nonlinear Schr\"{o}dinger equation in 2d with linear, cubic and 
quintic damping terms. Our numerical results show that quintic 
damping always arrests blowup, whereas linear and cubic damping 
can arrest blowup only when the damping parameter $\dt$ is bigger 
than a certain threshold value $\dt_{\rm th}$. We will apply this 
novel method to solve the 3d Gross-Pitaevskii equation with a 
quintic damping term and compare the numerical results with the 
experimental dynamics \cite{Donley} of collapsing and exploding 
BECs \cite{BaoBEC}. 
 
\bigskip 
 
\begin{center} 
{\large \bf Acknowledgment} 
\end{center} 
W.B. acknowledges support  by the National University of 
Singapore grant No. R-151-000-027-112. 
This work was supported by the WITTGENSTEIN-AWARD of P. Markowich 
and P. Zoller which is funded by the Austrian National Science 
Foundation FWF. The authors also acknowledge hospitality of the 
International Erwin Schr\"{o}dinger Institute in Vienna where this 
work was initiated.

 


\newpage 
 
\begin{figure}[htb] 
\centerline{a) \psfig{figure=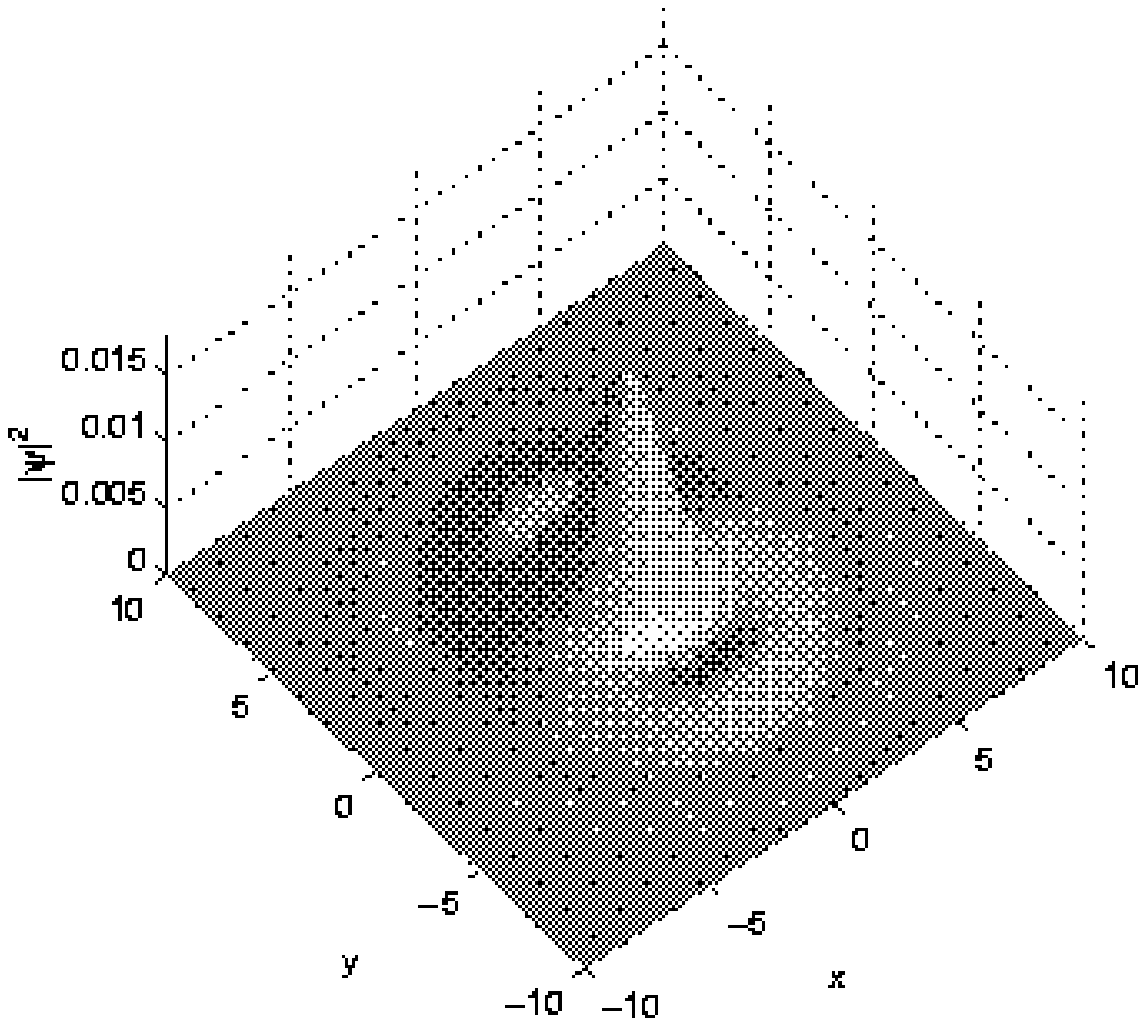,height=6cm,width=7cm,angle=0} 
\quad b) \psfig{figure=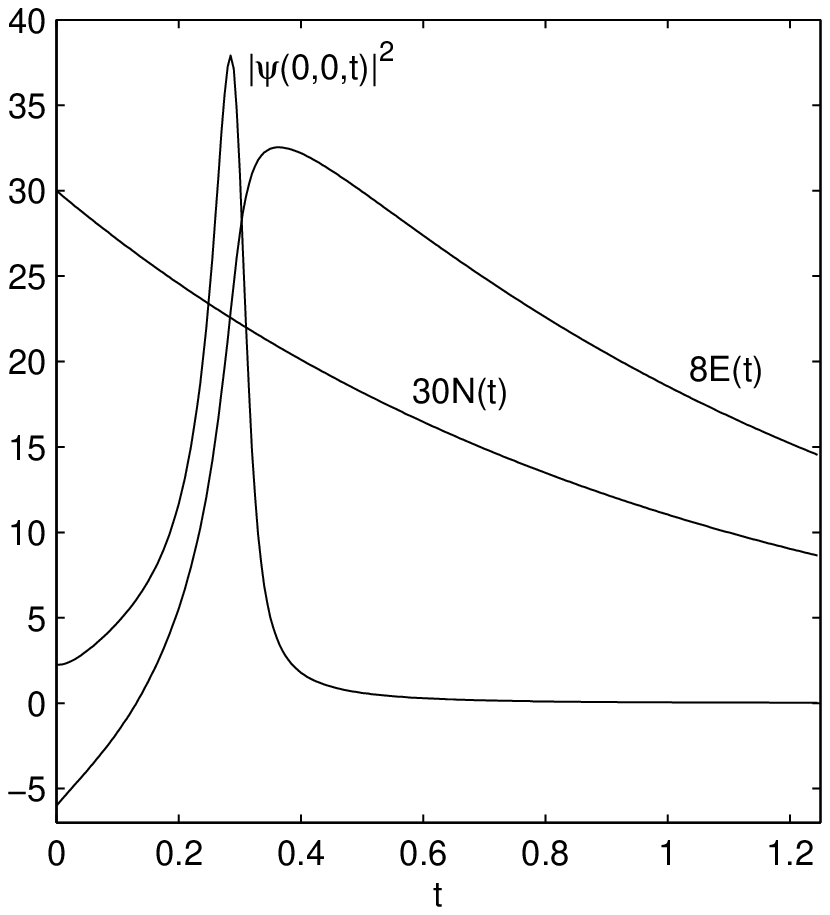,height=6cm,width=6cm,angle=0}} 
\centerline{c)\psfig{figure=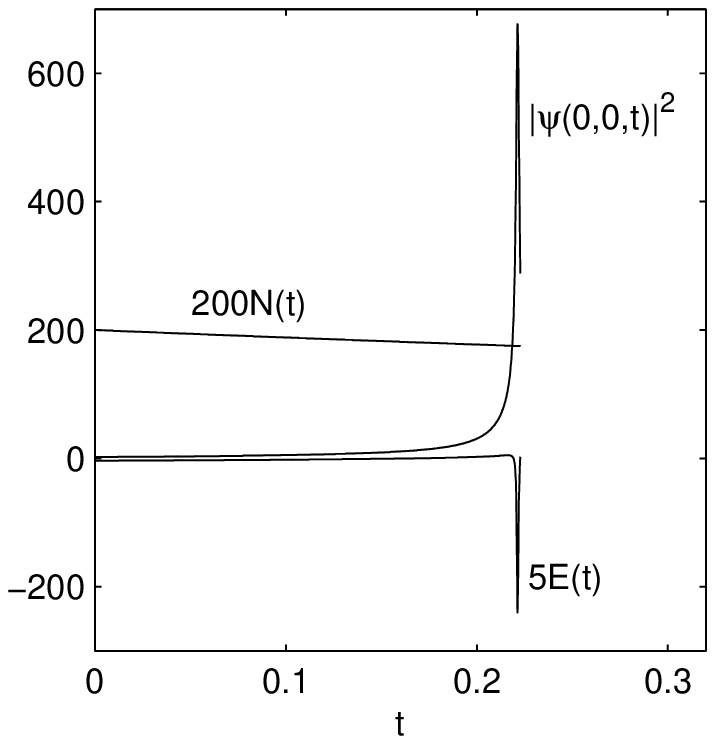,height=6cm,width=7cm,angle=0} 
\quad d)\psfig{figure=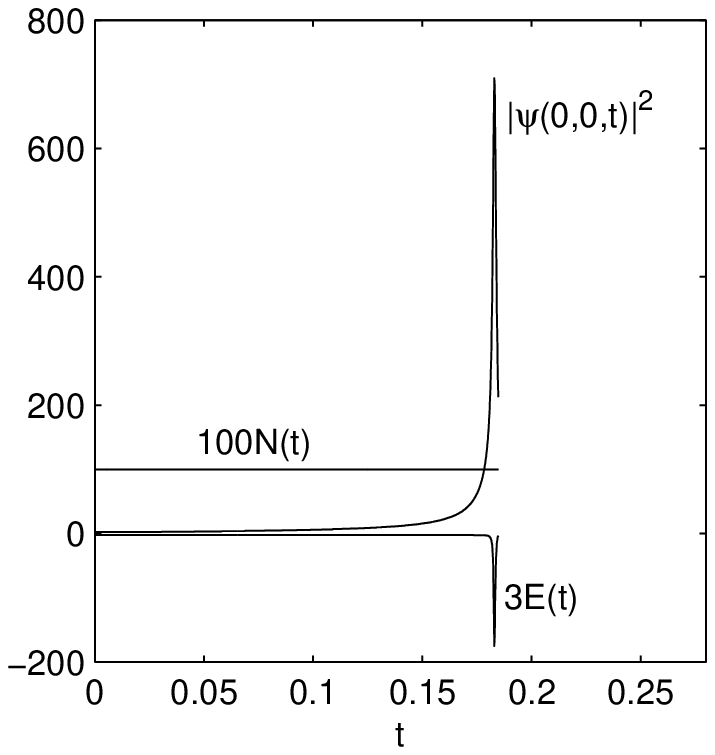,height=6cm,width=7cm,angle=0}} 
 \centerline{e)\psfig{figure=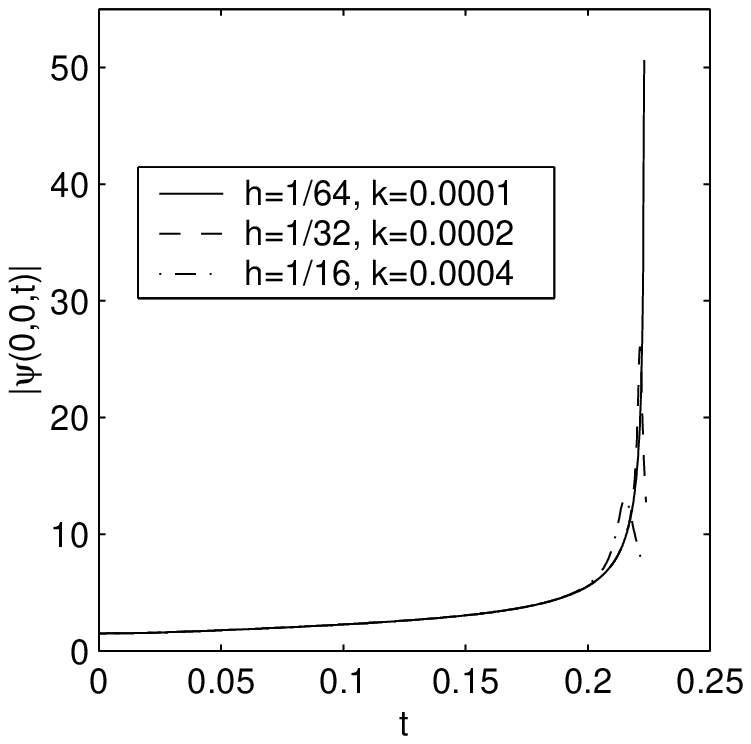,height=6cm,width=7cm,angle=0} 
\quad f)\psfig{figure=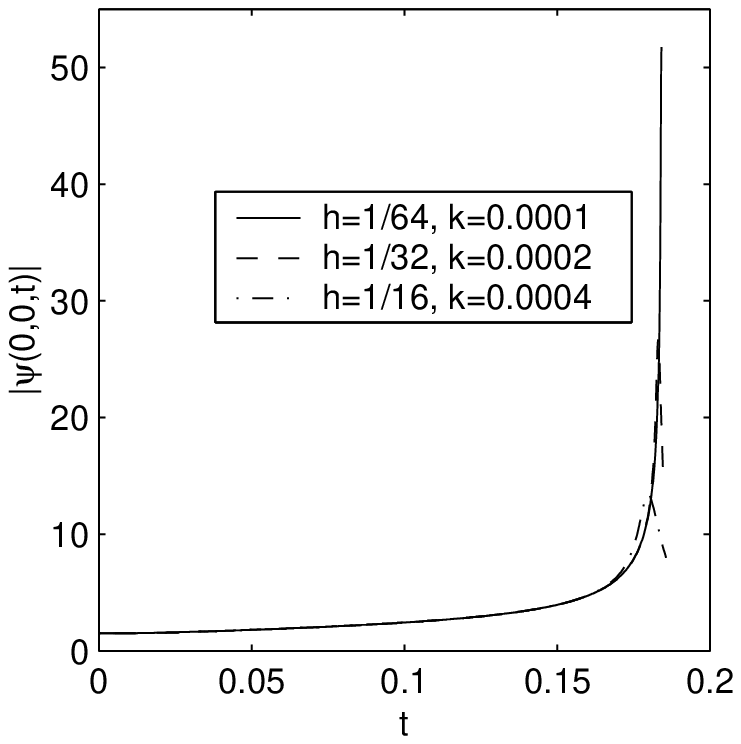,height=6cm,width=7cm,angle=0}}

Figure 1: Numerical results  in Example 1 case I. \quad 
a). Surface plot of the density $|\psi|^2$ at 
time $t=1.25$  with $\dt=0.5$. \quad 
Normalization, energy and central density 
$|\psi(0,0,t)|^2$ as functions of time: b). with 
$\dt=0.5$, c).  $\dt=0.3$, 
d).  $\dt=0$ (no damping). 
\quad Blowup study: e).  $\dt=0.3$, 
f).  $\dt=0$ (no damping). 
\end{figure} 
 
\begin{figure}[htb] 
\centerline{a).\psfig{figure=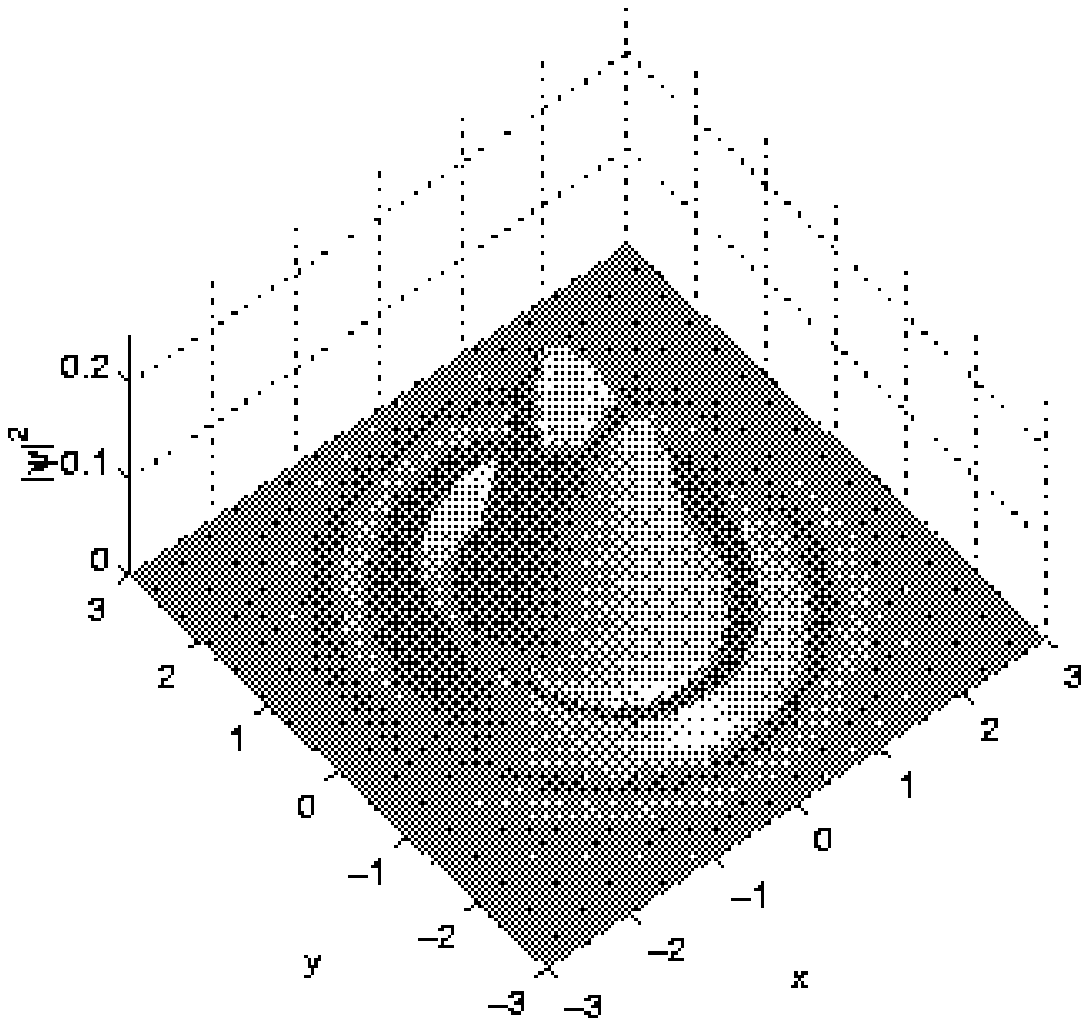,height=7cm,width=7cm,angle=0} 
\quad  b).\psfig{figure=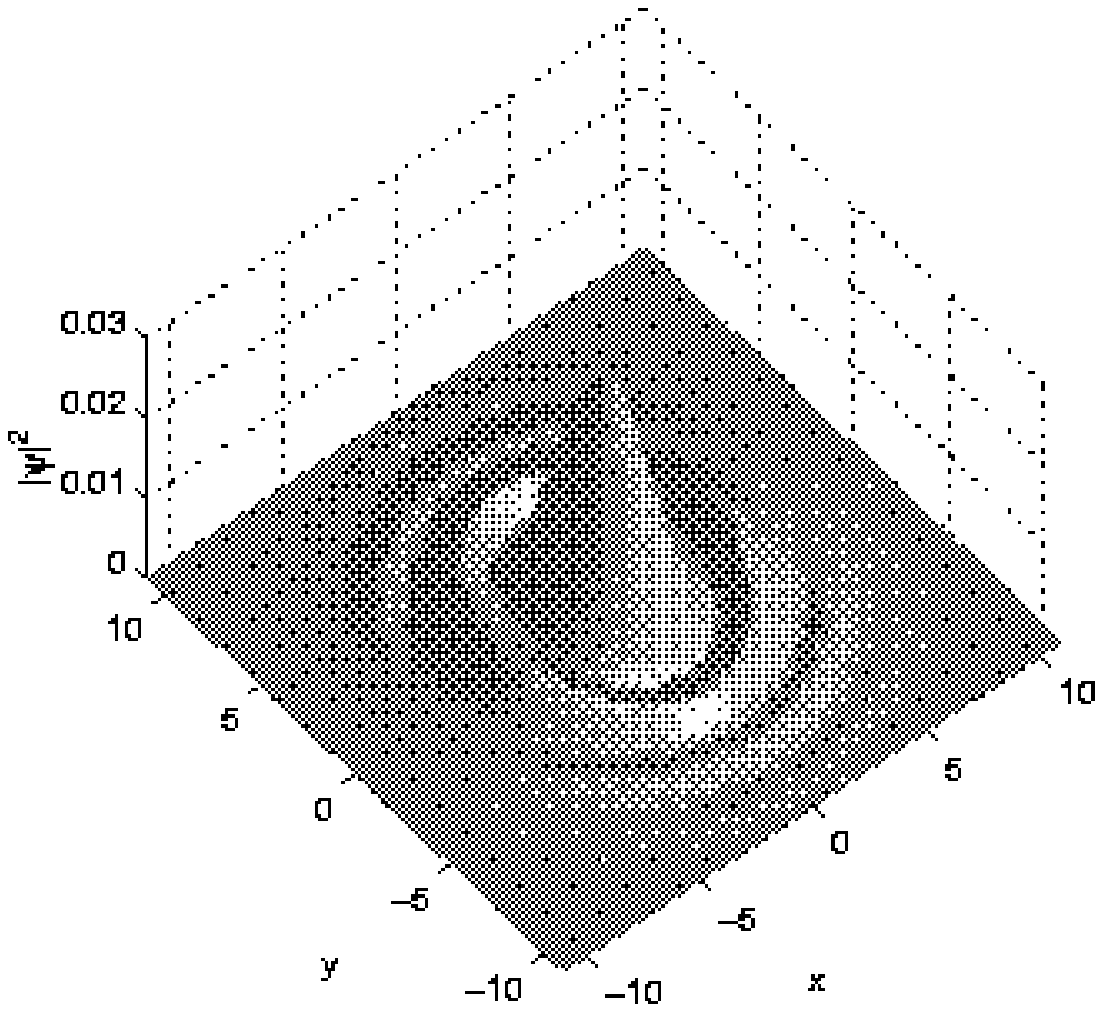,height=7cm,width=7cm,angle=0}} 
\bigskip 
\bigskip 
\bigskip 
\bigskip 
\centerline{c).\psfig{figure=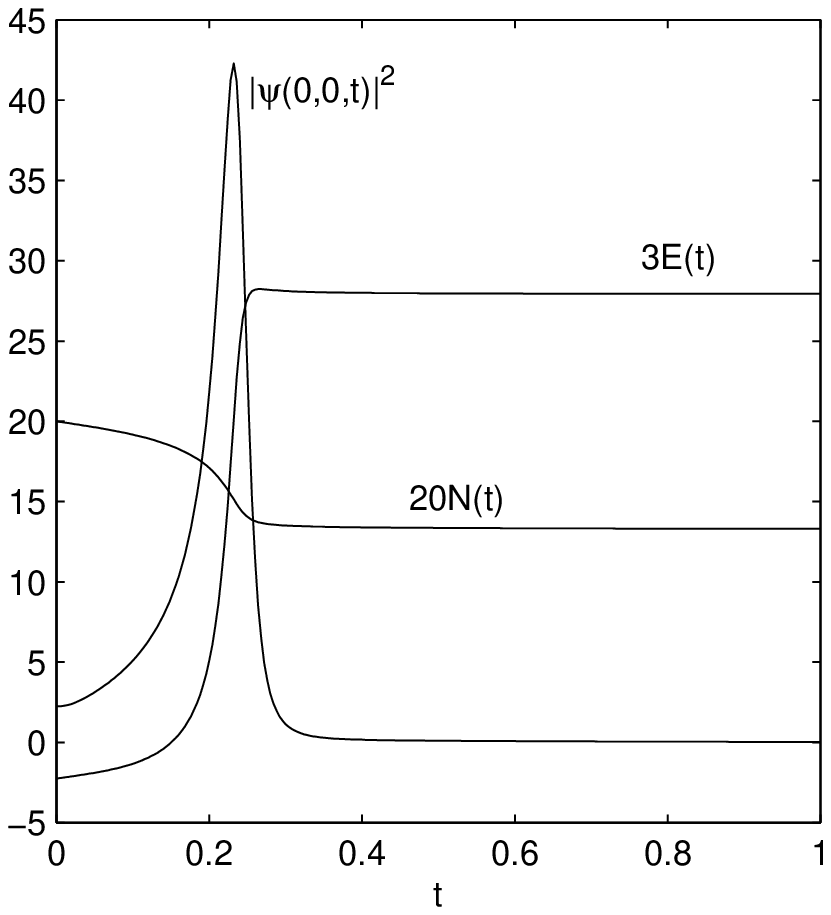,height=7cm,width=7cm,angle=0} 
\quad d).\psfig{figure=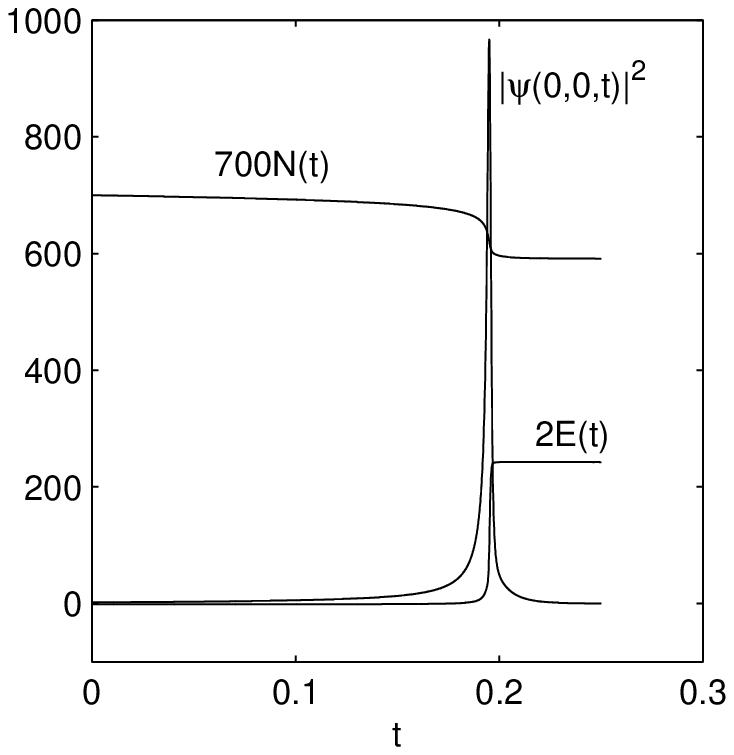,height=7cm,width=7cm,angle=0}} 
 
Figure 2:  Numerical results  in Example 1 case II. \quad 
 Surface plot of the density $|\psi|^2$ 
 with $\dt=0.02$: a). At time $t=0.4$, b). $t=1.0$. \quad 
Normalization, energy and central density 
$|\psi(0,0,t)|^2$ as functions of time: 
c).  with $\dt=0.02$, d). $\dt=0.005$ (with $h=1/128$, $k=0.00002$). 
 
\end{figure} 
 
\begin{figure}[htb] 
\centerline{a).\psfig{figure=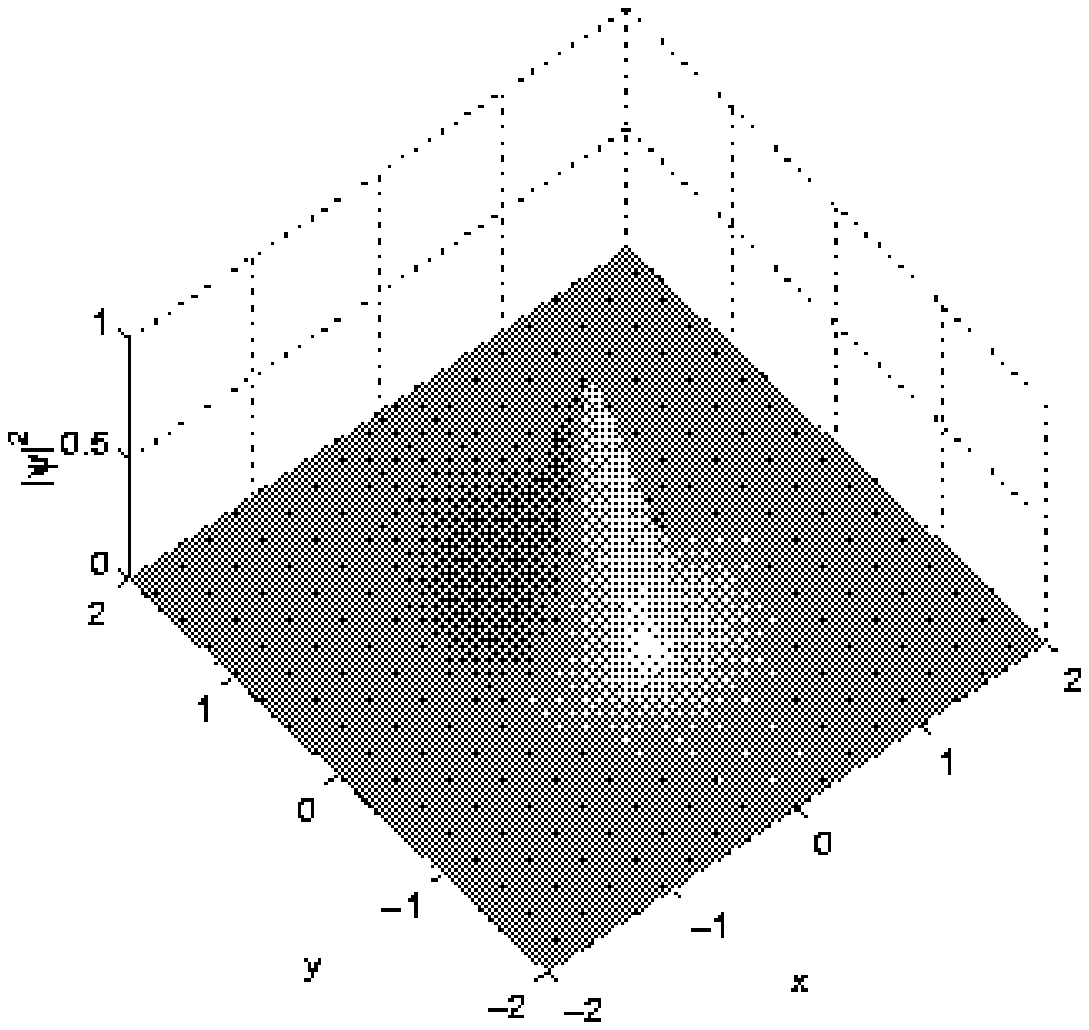,height=7cm,width=7cm,angle=0} 
\quad b).\psfig{figure=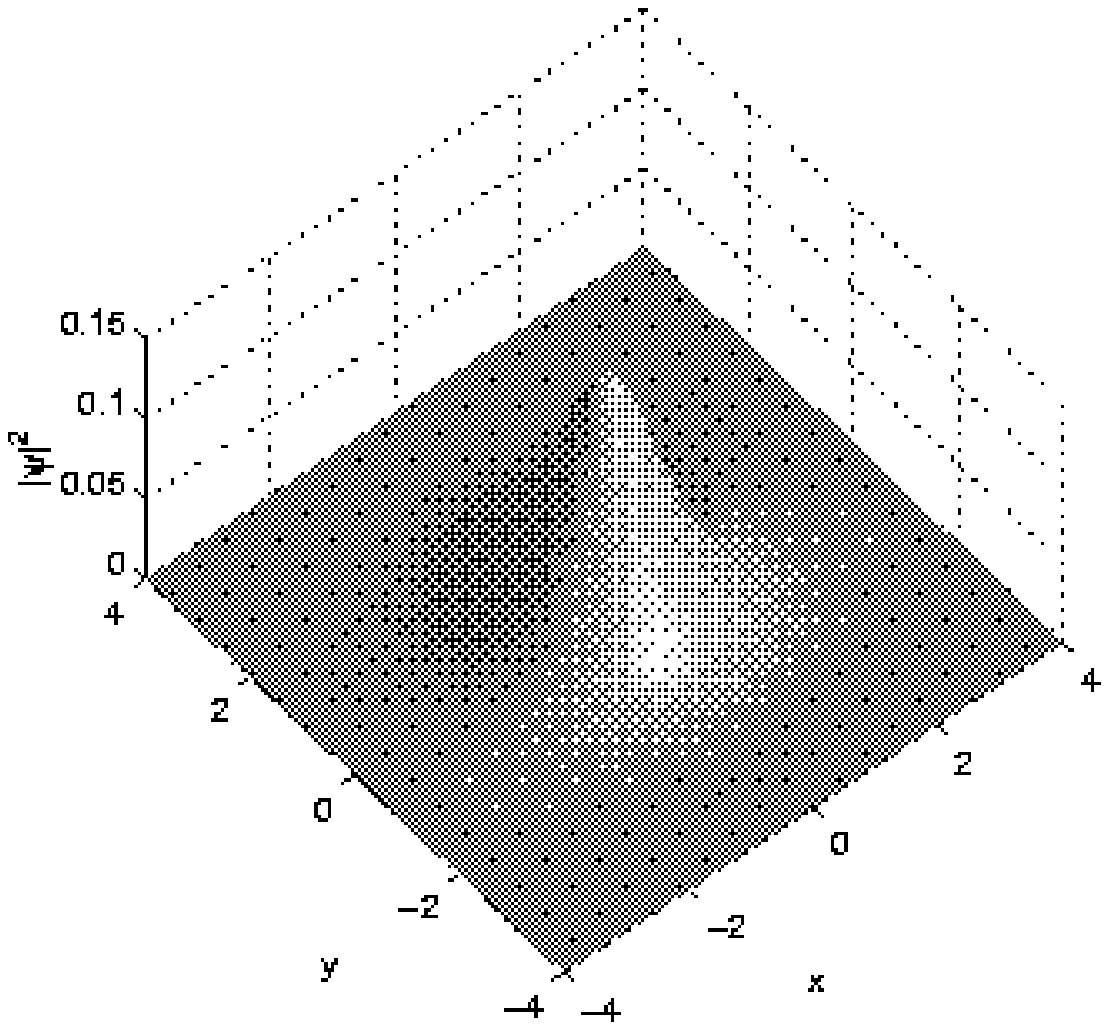,height=7cm,width=7cm,angle=0}} 
\bigskip 
\bigskip 
\bigskip 
\bigskip 
\centerline{c).\psfig{figure=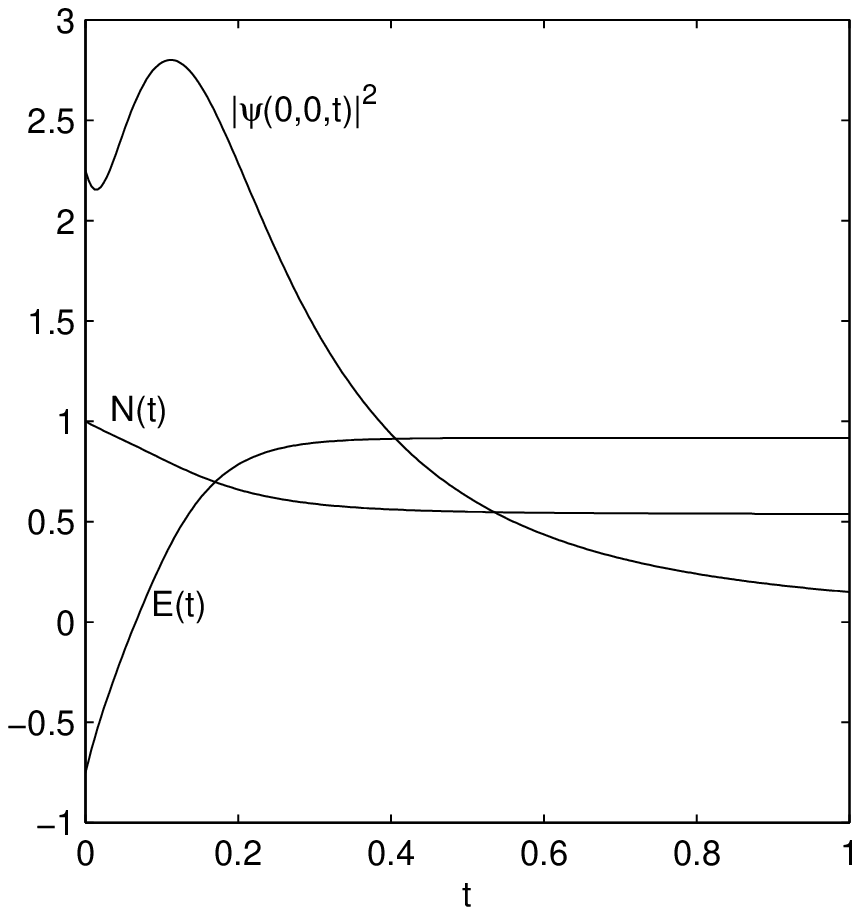,height=7cm,width=7cm,angle=0} 
\quad d).\psfig{figure=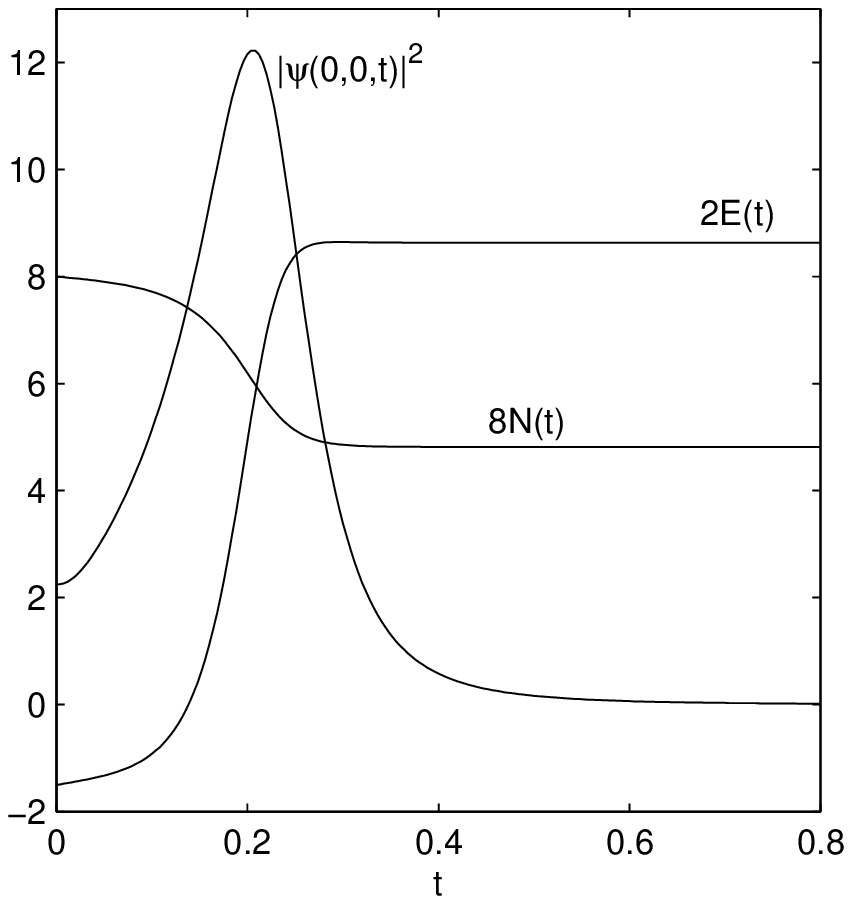,height=7cm,width=7cm,angle=0}} 
 
Figure 3:  Numerical results  in Example 1 case III. \quad 
 Surface plot of the density $|\psi|^2$ 
 with $\dt=0.01$: a). At time $t=0.4$, b). $t=1.0$. \quad 
Normalization, energy and central density 
$|\psi(0,0,t)|^2$ as functions of time: 
c).  with $\dt=0.01$, d). $\dt=0.001$. 
 
\end{figure}

\begin{figure}[htb] 
\centerline{a).\psfig{figure=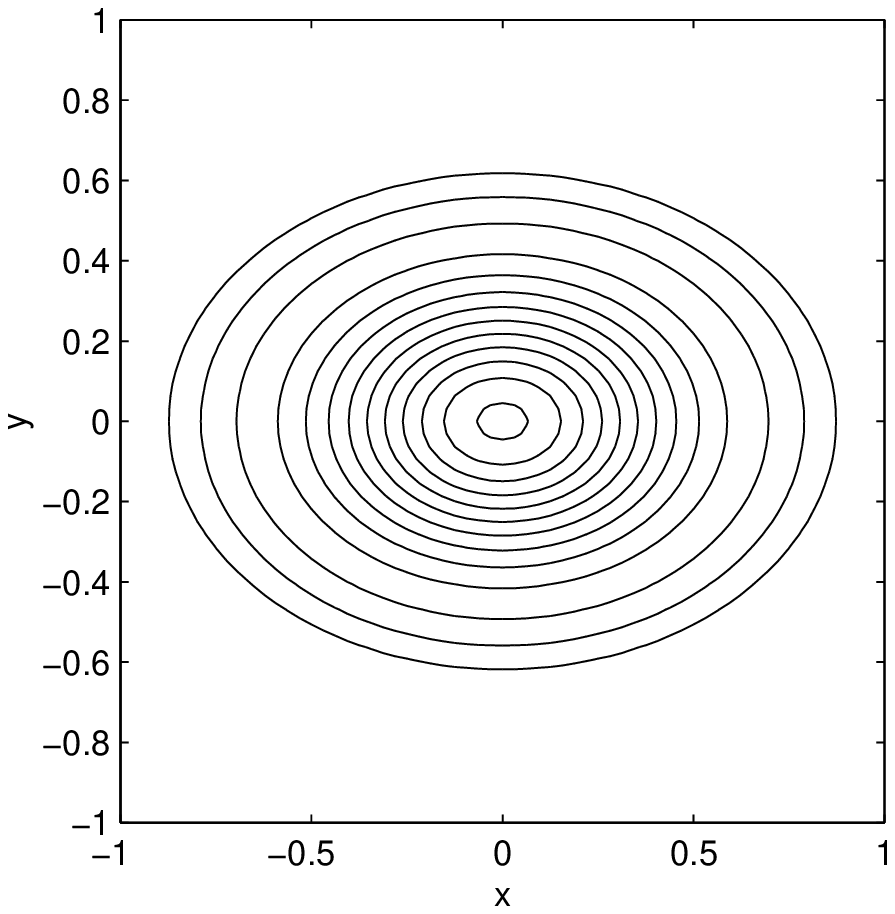,height=6cm,width=7cm,angle=0} 
\quad b)\psfig{figure=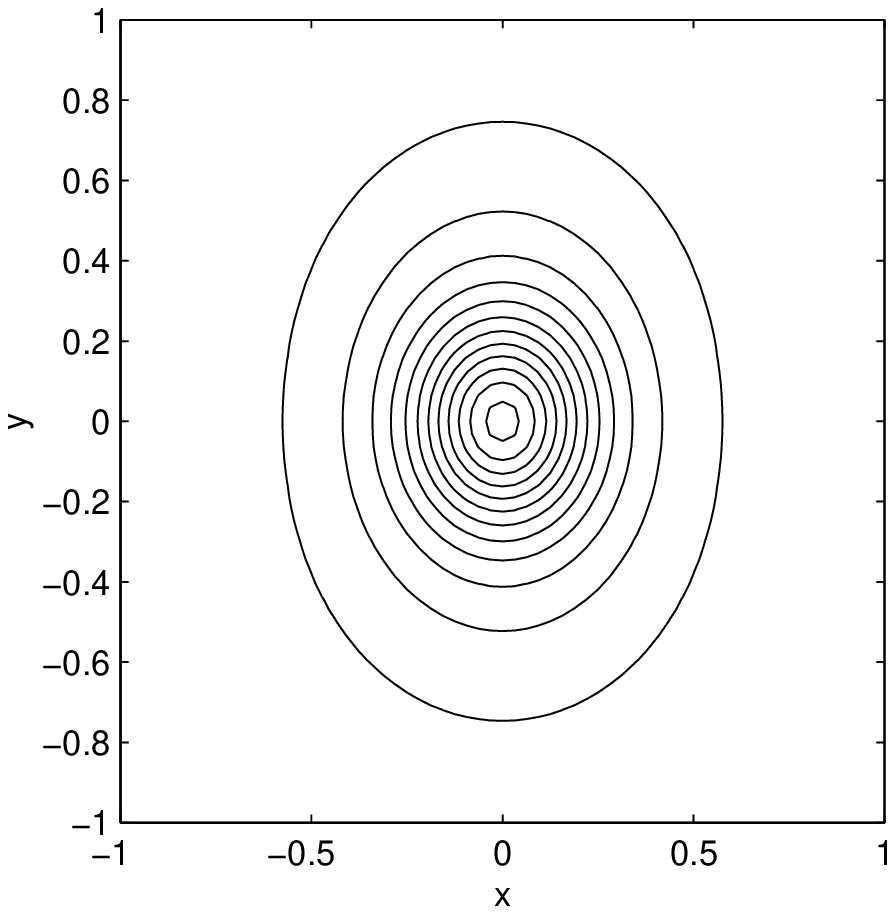,height=6cm,width=7cm,angle=0} } 
\centerline{c).\psfig{figure=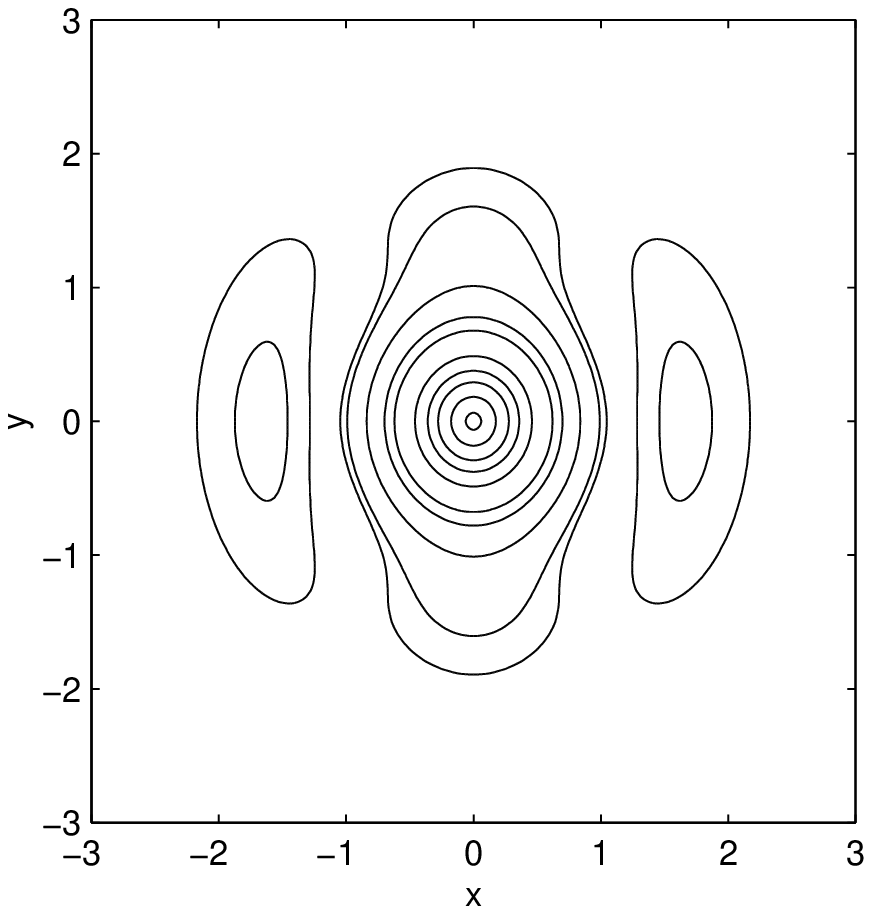,height=6cm,width=7cm,angle=0} 
\quad d)\psfig{figure=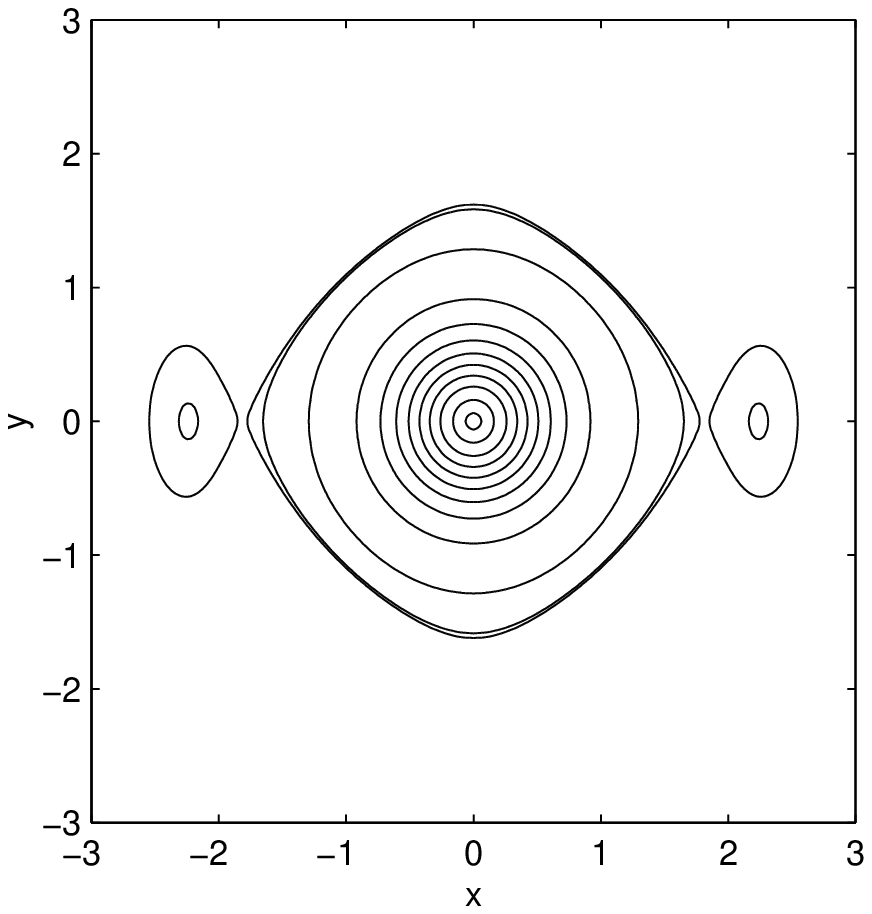,height=6cm,width=7cm,angle=0} } 
\centerline{e).\psfig{figure=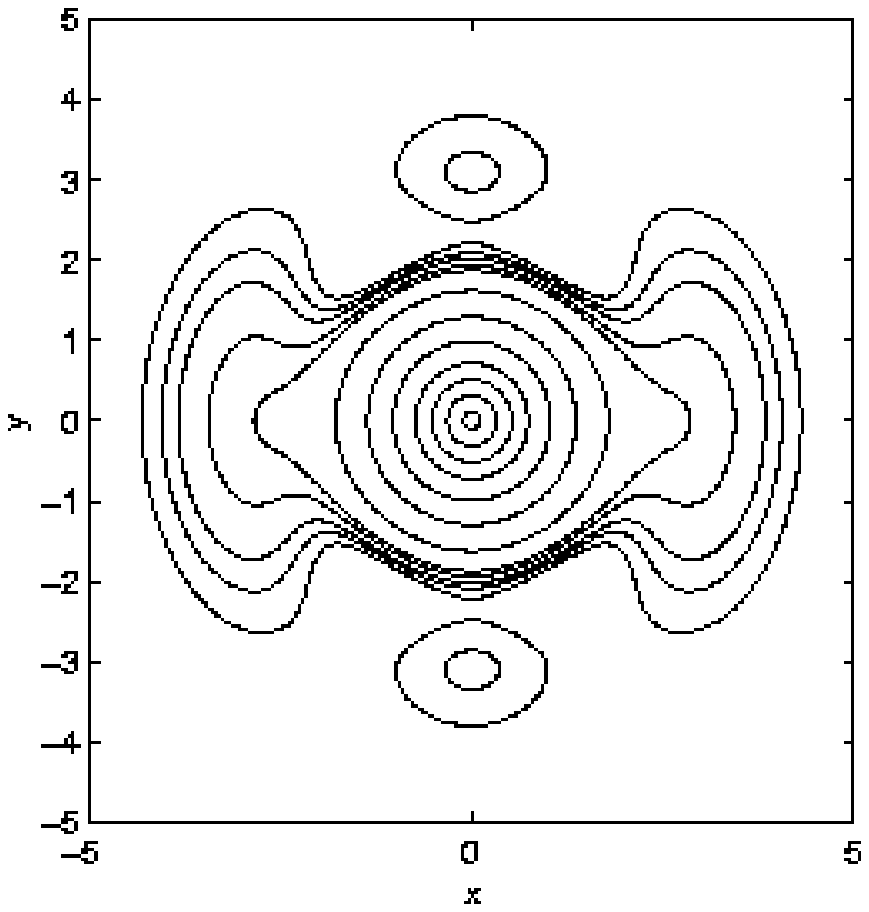,height=6cm,width=7cm,angle=0} 
\quad f)\psfig{figure=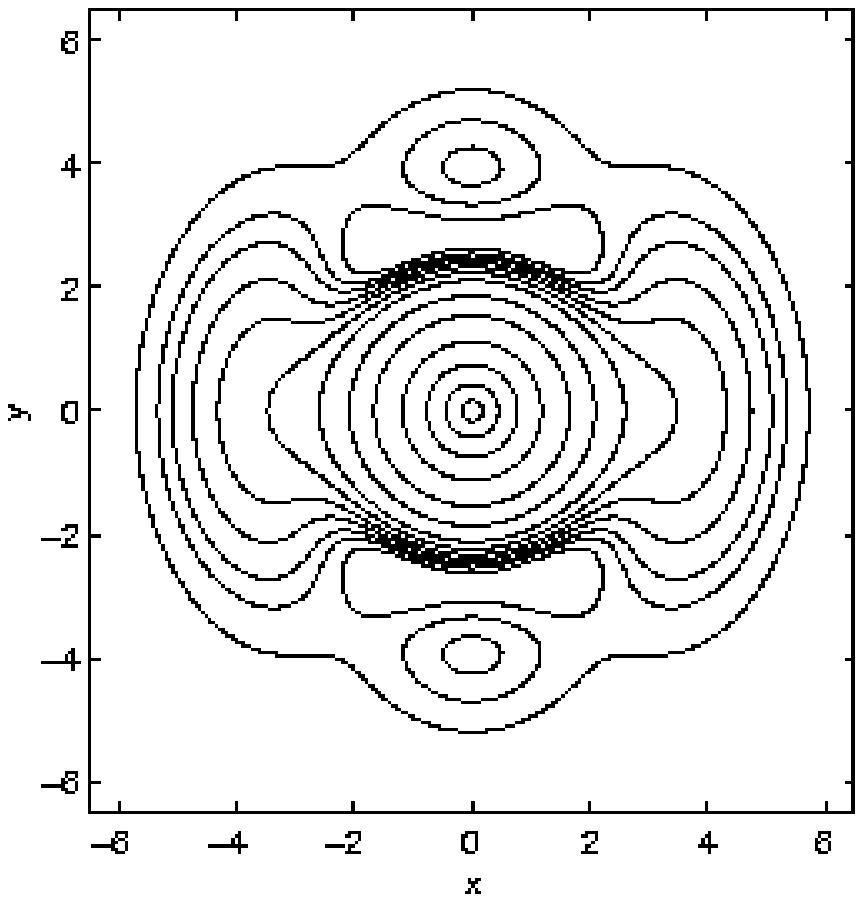,height=6cm,width=7cm,angle=0} } 
 
 Figure 4:  Contour plots of the density $|\psi|^2$ 
at different times in Example 1 case III with $\dt=0.01$. 
\  a). $t=0$, 
b). $t=0.2$, c). $t=0.4$, d). $t=0.6$, e). $t=0.8$, f). $t=1$. 
\end{figure}

\begin{figure}[htb] 
\centerline{a).\psfig{figure=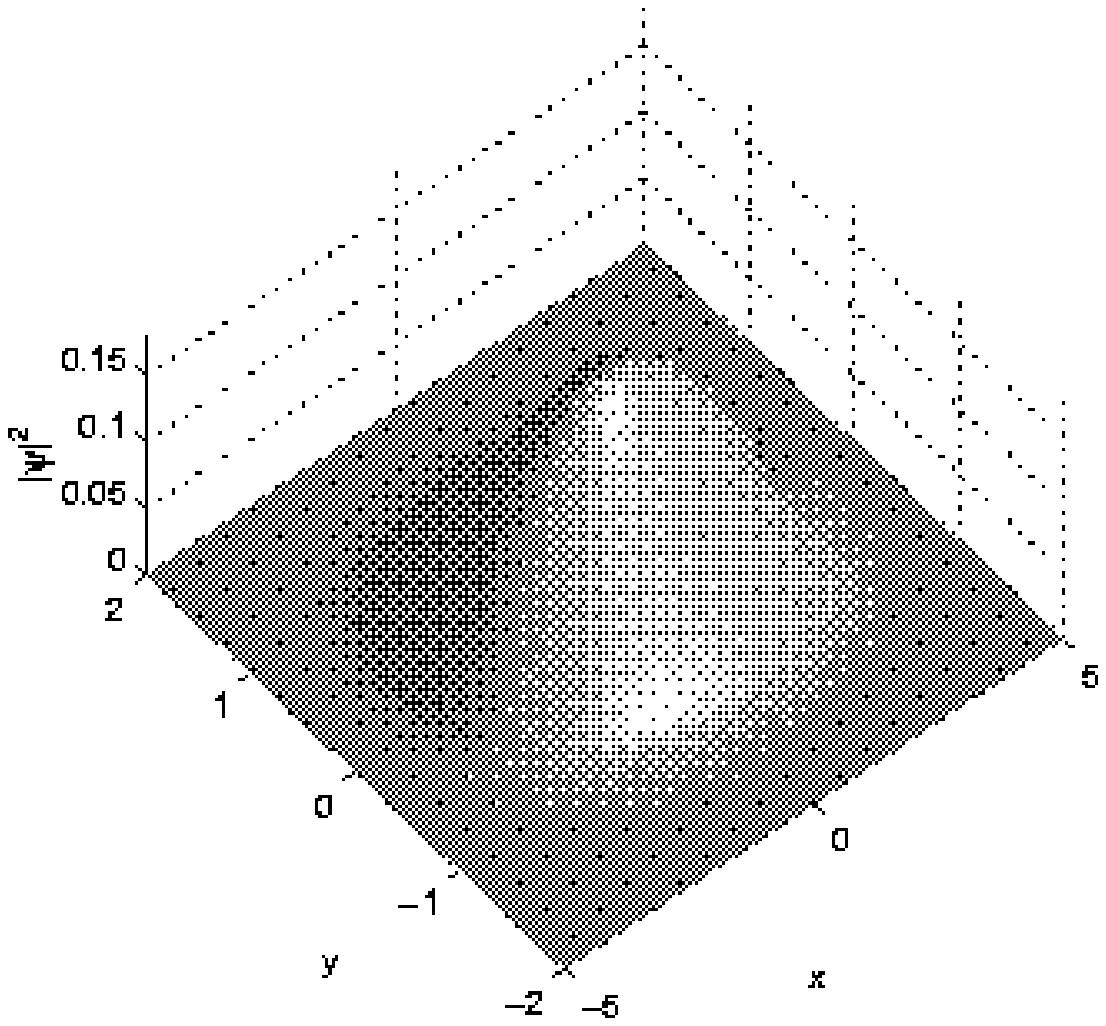,height=6cm,width=7cm,angle=0} 
\quad b).\psfig{figure=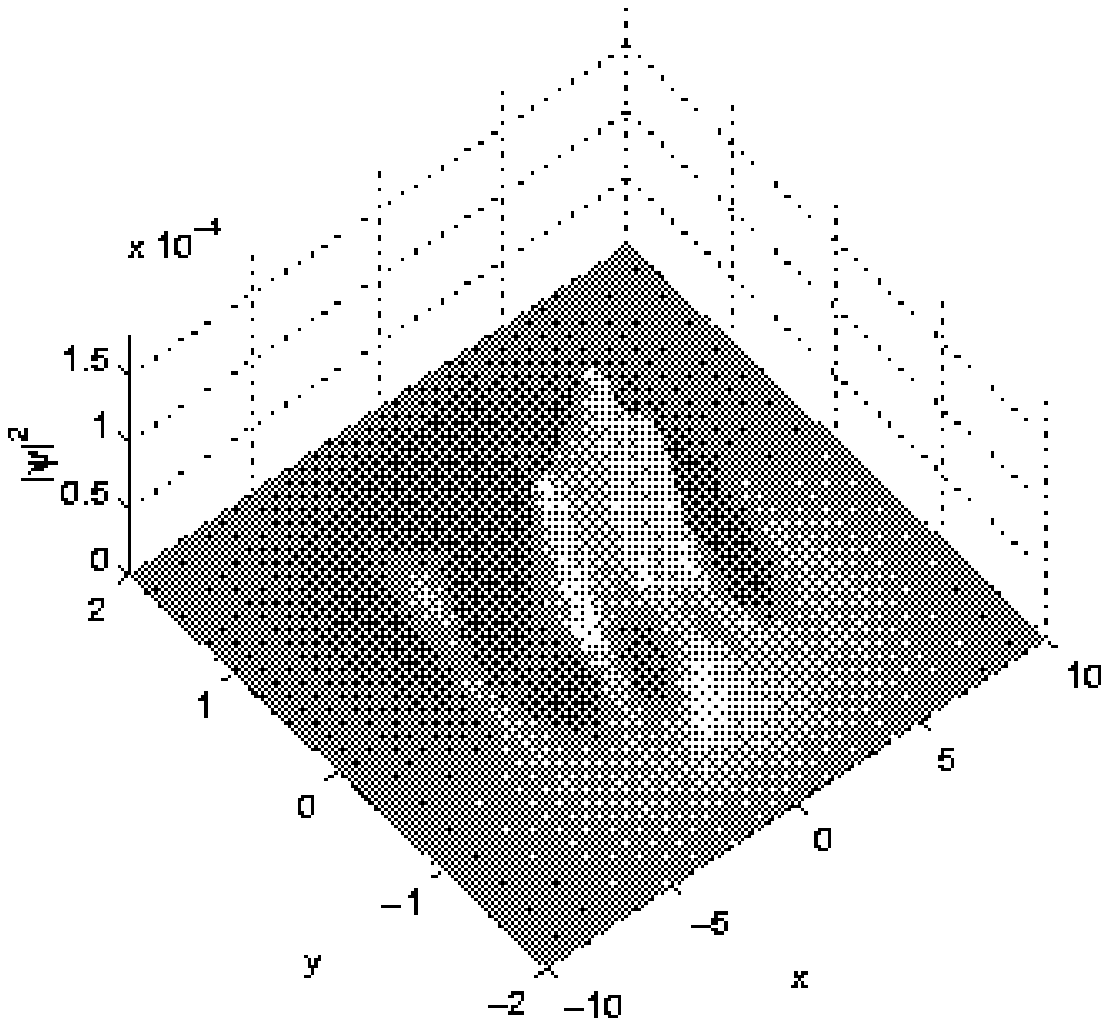,height=6cm,width=7cm,angle=0}} 
\centerline{c).\psfig{figure=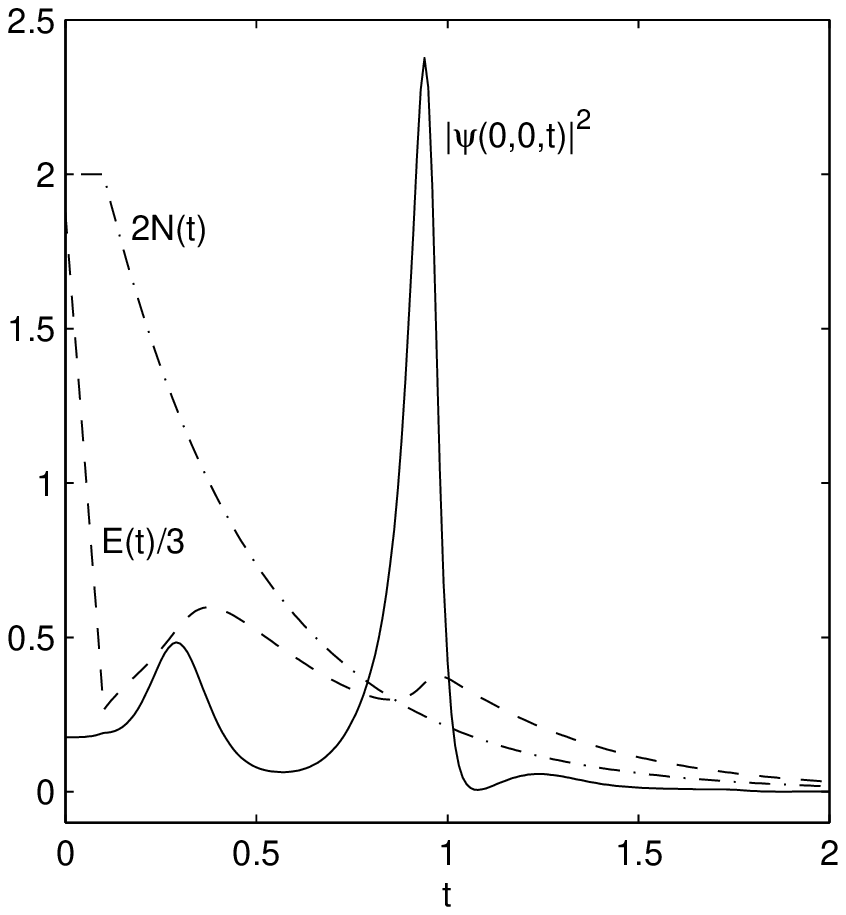,height=6cm,width=7cm,angle=0} 
\quad d)\psfig{figure=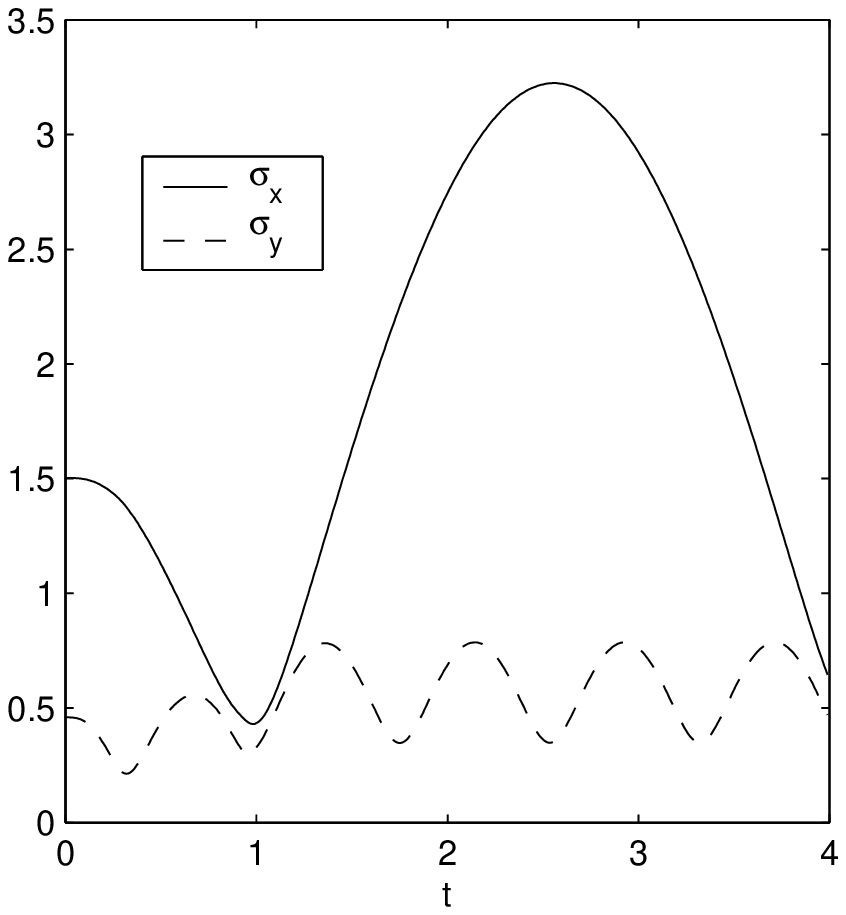,height=6cm,width=7cm,angle=0} } 
\centerline{e).\psfig{figure=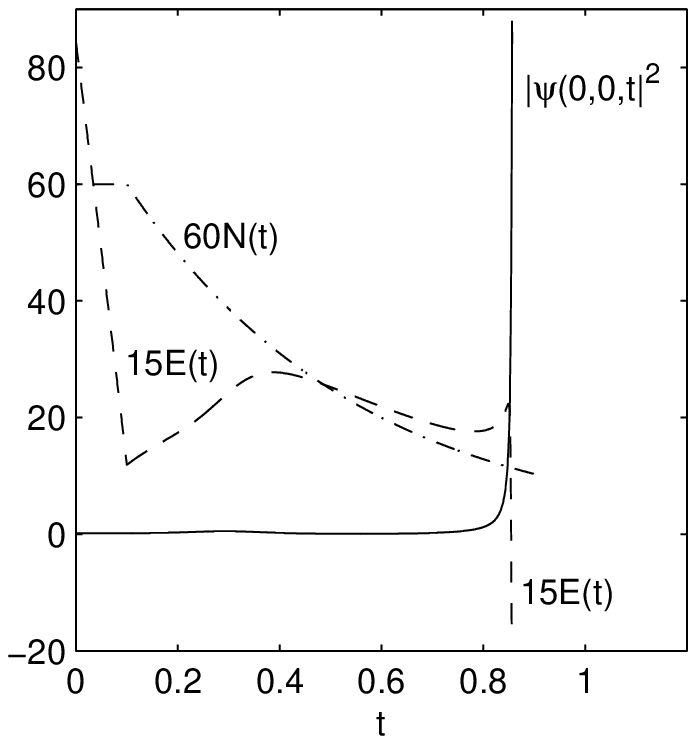,height=6cm,width=7cm,angle=0} 
\quad f)\psfig{figure=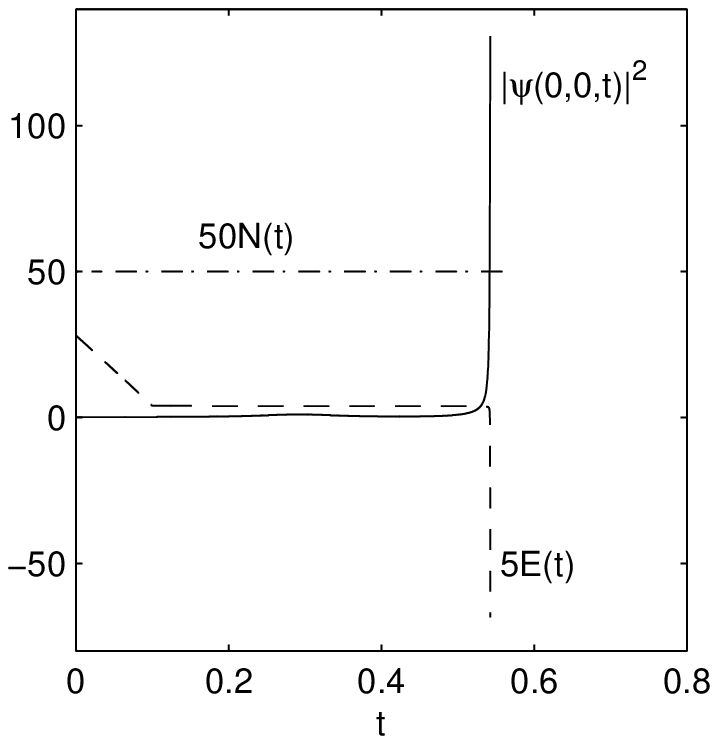,height=6cm,width=7cm,angle=0} } 
 
Figure 5:  Numerical results  in Example 2 case I. \quad 
 Surface plot of the density $|\psi|^2$ 
 with $\dt=1.25$: a). At time $t=0$ (ground-state solution), 
b). $t=2.8$.  \quad 
Normalization, energy and central density 
$|\psi(0,0,t)|^2$ as functions of time: 
c).  with $\dt=1.25$, e). $\dt=1.1$, f). $\dt=0$ (no damping). \quad 
d). Condensate widths with $\dt=1.25$. 
 
\end{figure} 
 
\begin{figure}[htb] 
\centerline{a).\psfig{figure=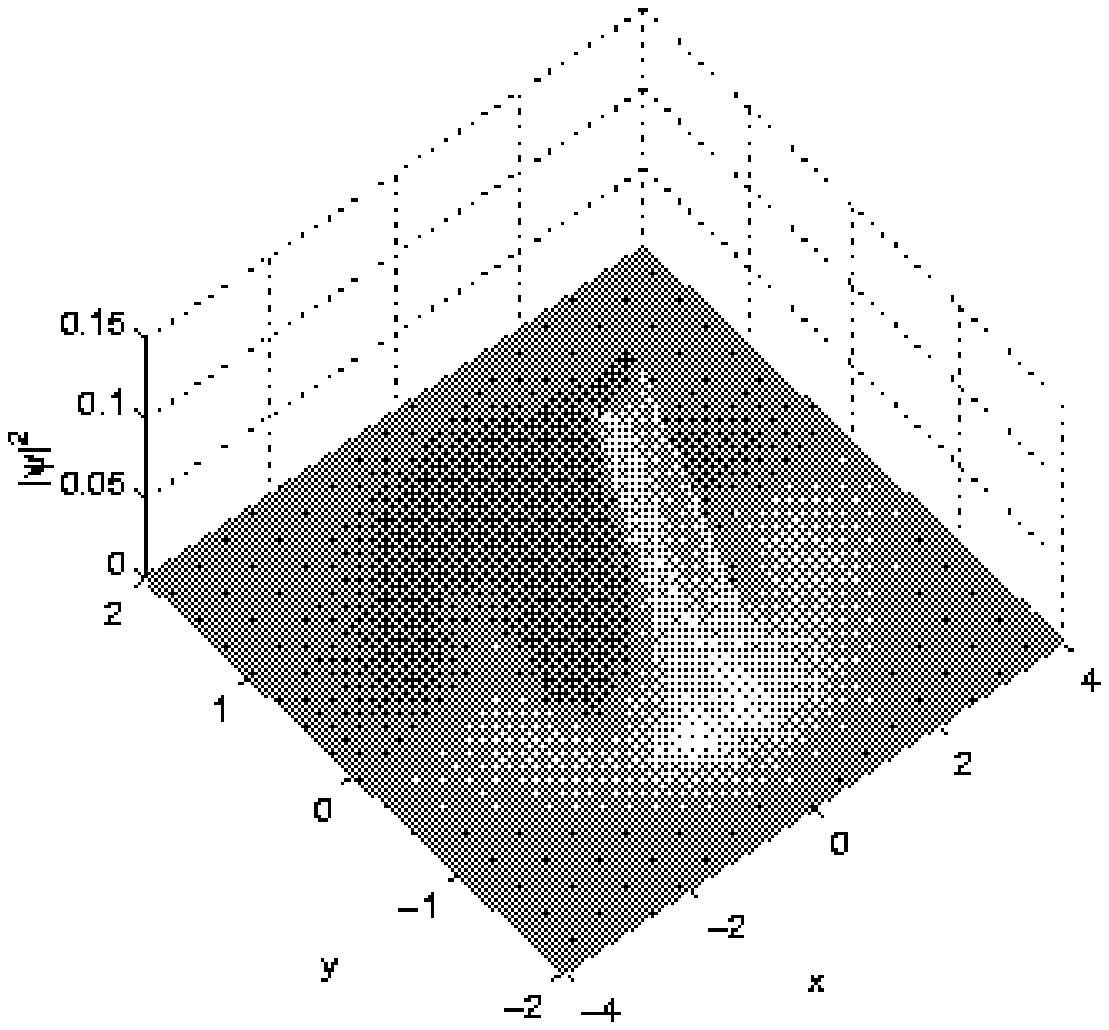,height=6cm,width=7cm,angle=0} 
\quad \psfig{figure=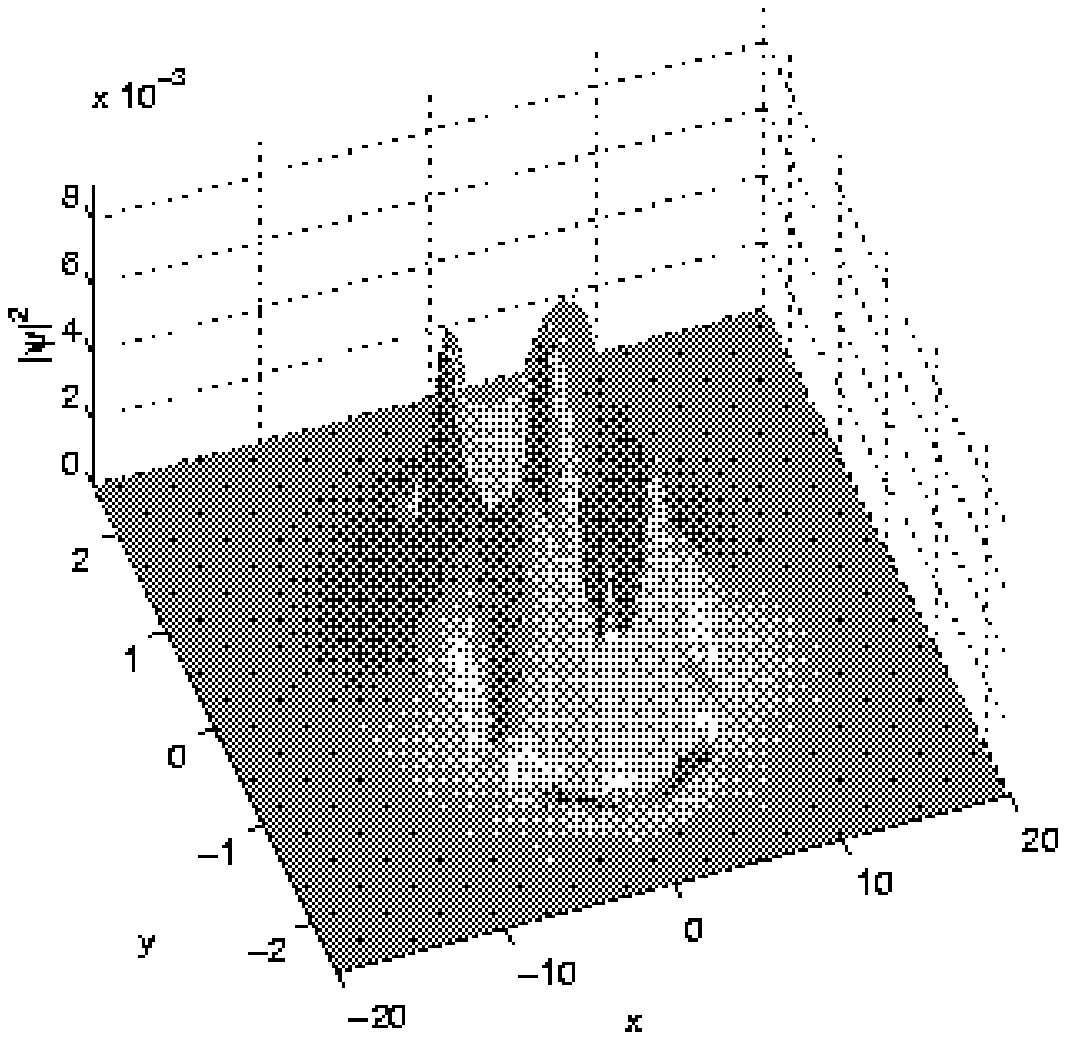,height=6cm,width=7cm,angle=0}} 
\centerline{b).\psfig{figure=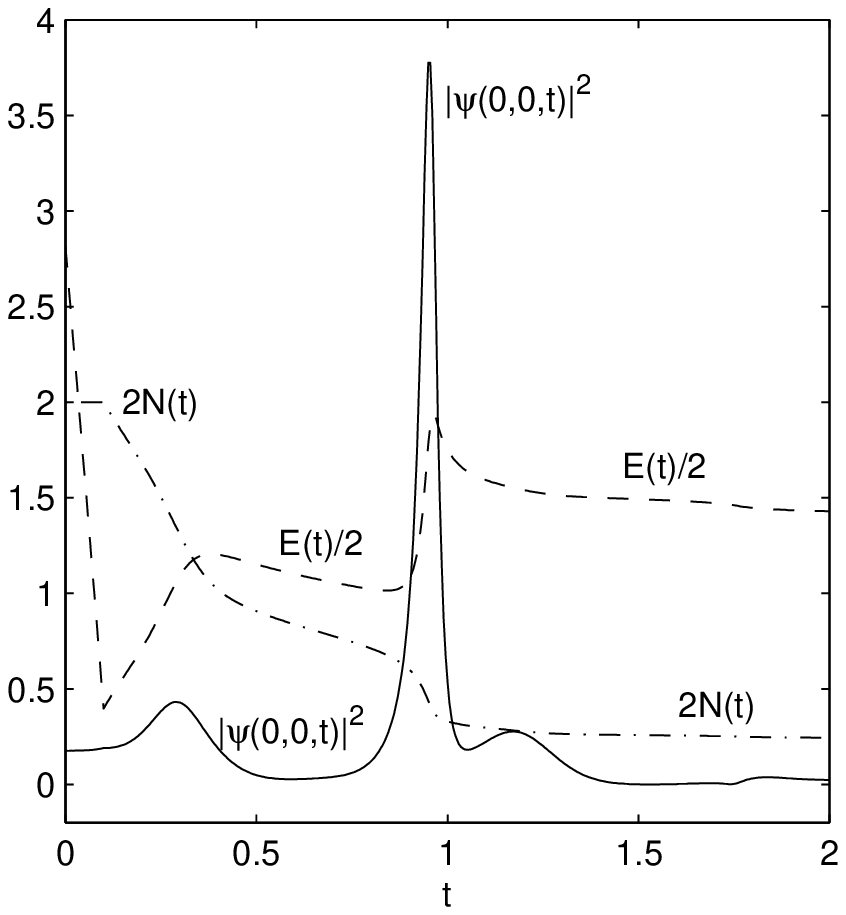,height=6cm,width=7cm,angle=0} 
\quad \psfig{figure=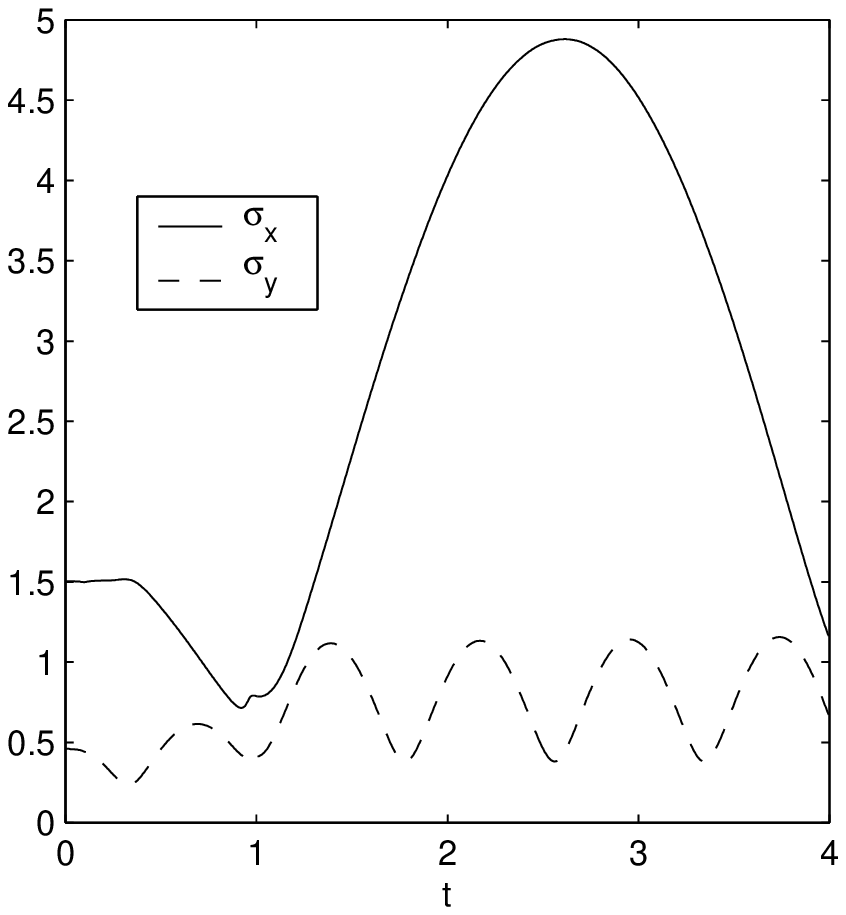,height=6cm,width=7cm,angle=0} } 
\centerline{c).\psfig{figure=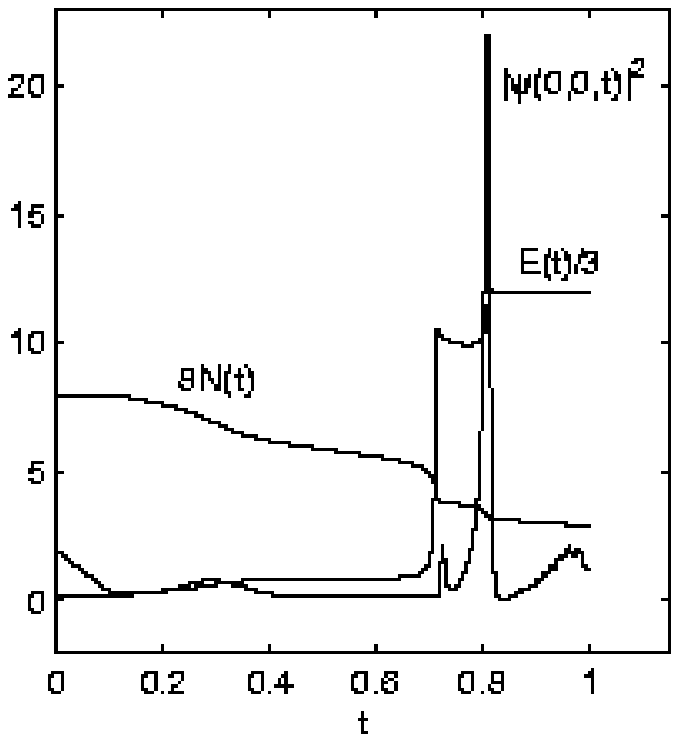,height=6cm,width=7cm,angle=0} 
\quad \psfig{figure=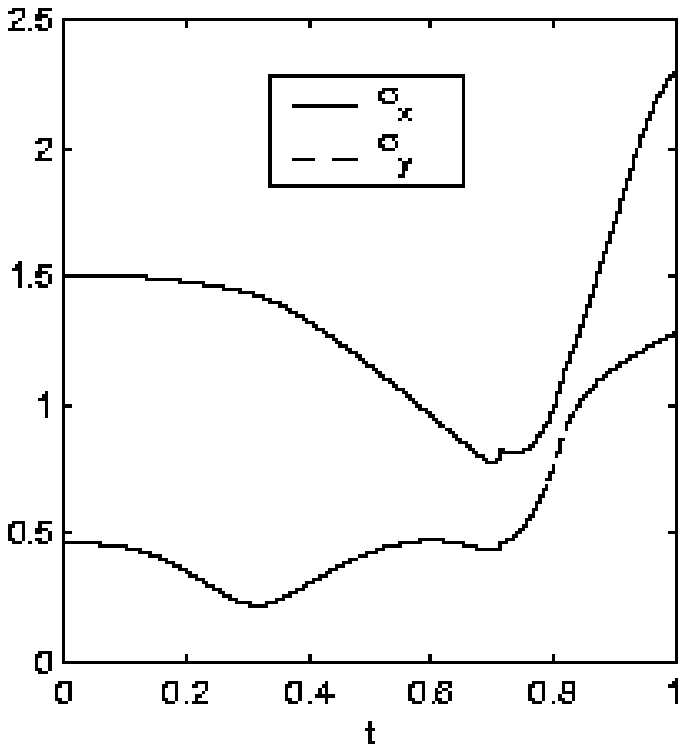,height=6cm,width=7cm,angle=0} } 
 
Figure 6:  Numerical results  in Example 2 case II. \quad 
 a). Surface plot of the density $|\psi|^2$ 
 with $\dt=0.15$: At time $t=0.8$ (left column) and 
$t=2.4$ (right column). \quad 
 Normalization, energy and central density 
 $|\psi(0,0,t)|^2$ (left column) 
and condensate widths (right column) 
as functions of time:  b).  With $\dt=0.15$;  c). $\dt=0.04$ 
(with . 
 
\end{figure} 
 
\begin{figure}[htb] 
\centerline{a).\psfig{figure=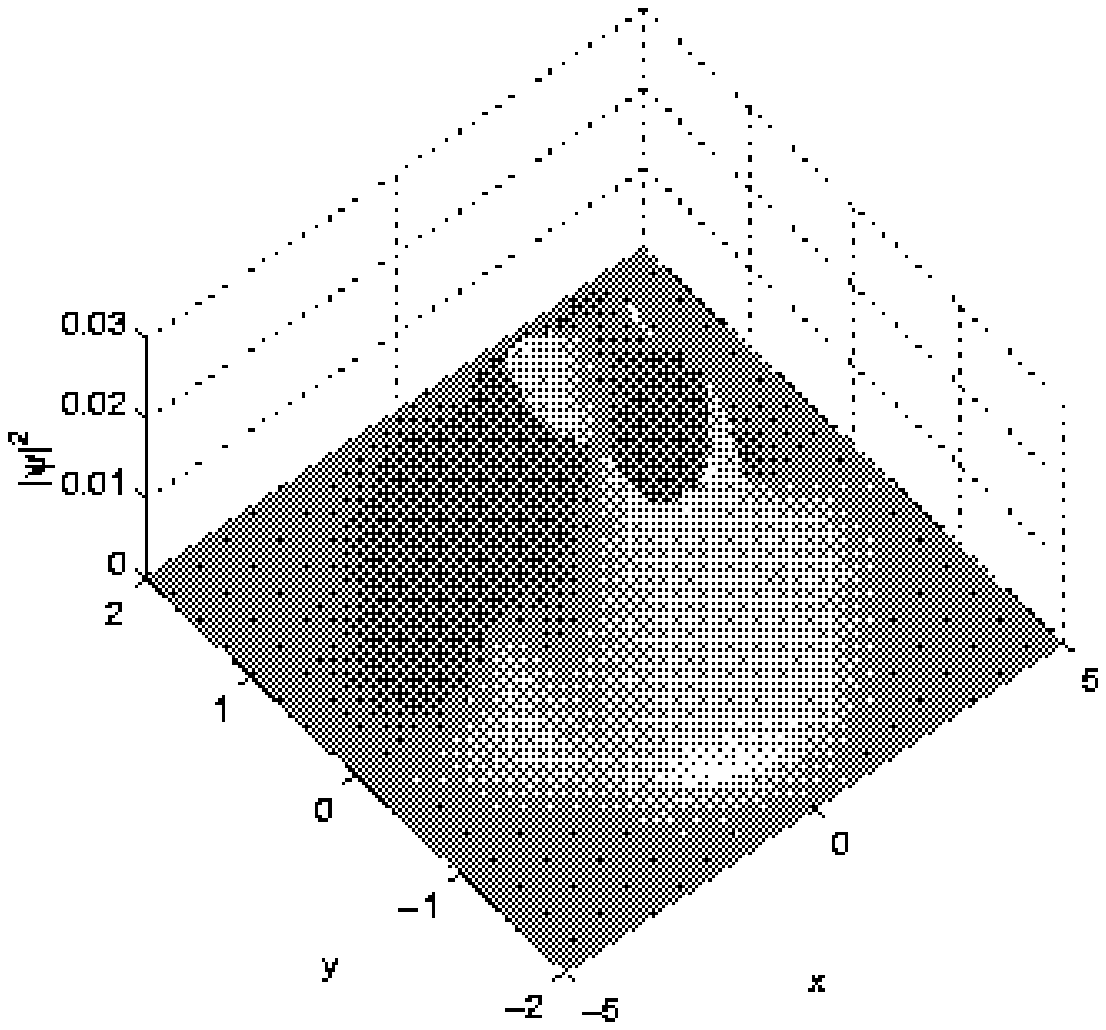,height=6cm,width=7cm,angle=0} 
\quad \psfig{figure=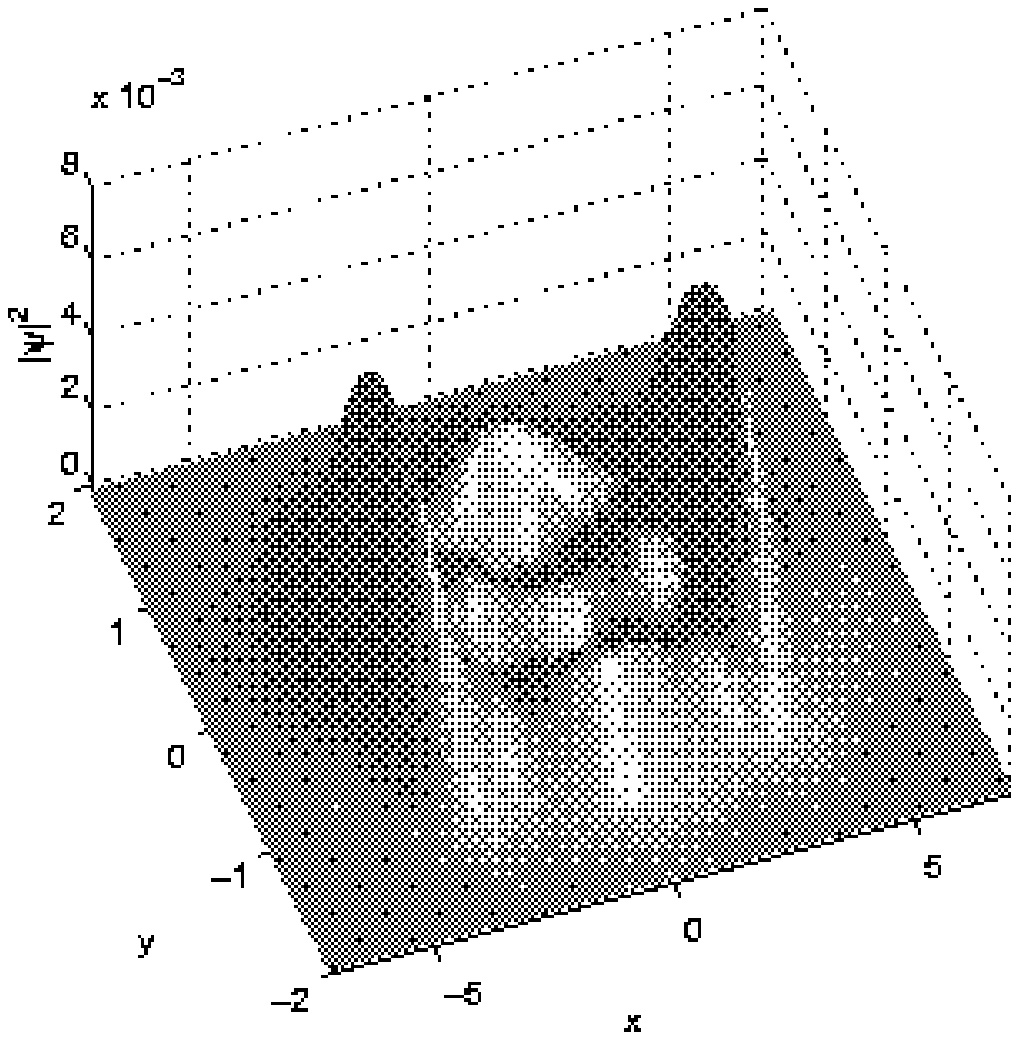,height=6cm,width=7cm,angle=0}} 
\centerline{b).\psfig{figure=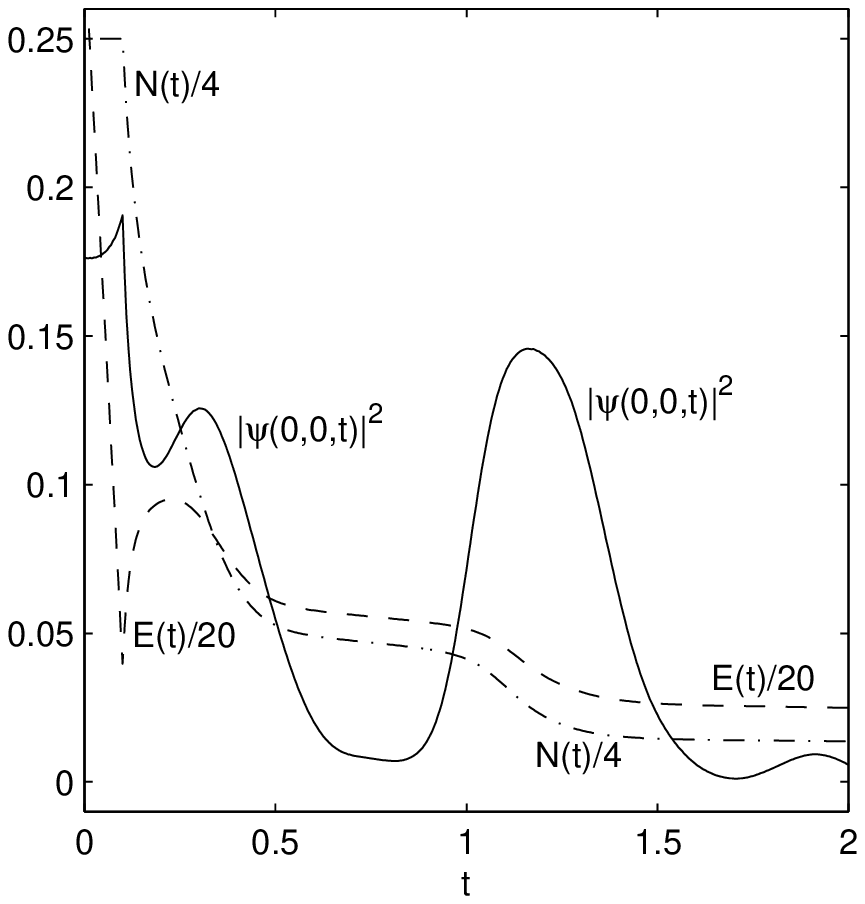,height=6cm,width=7cm,angle=0} 
\quad \psfig{figure=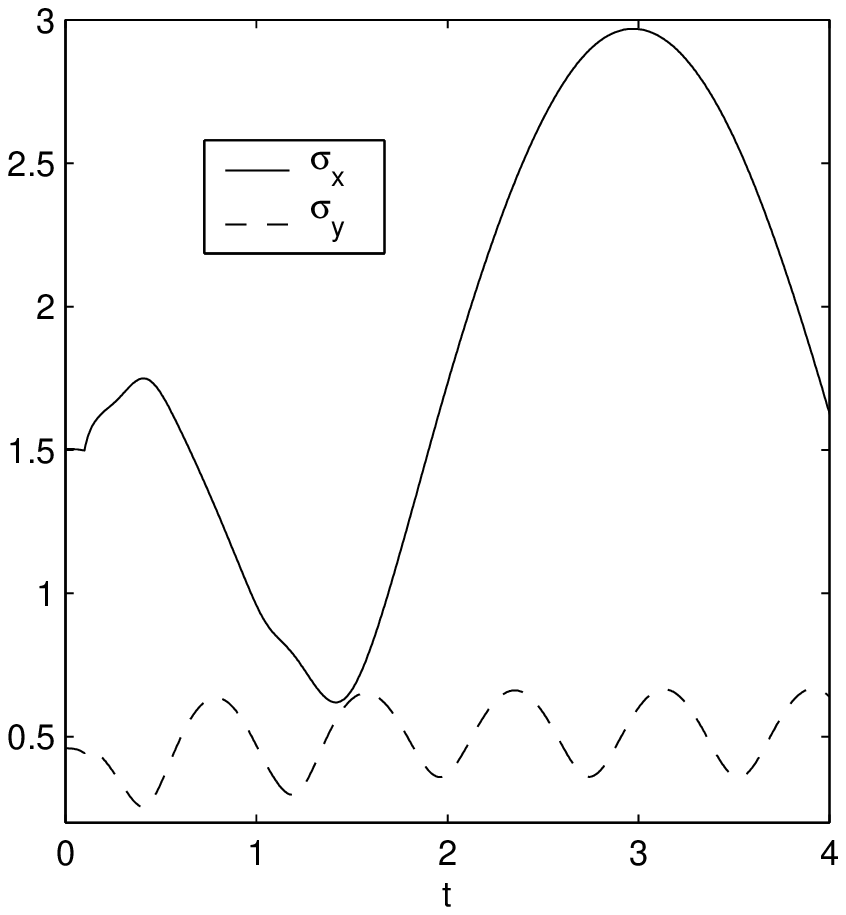,height=6cm,width=7cm,angle=0} } 
\centerline{c).\psfig{figure=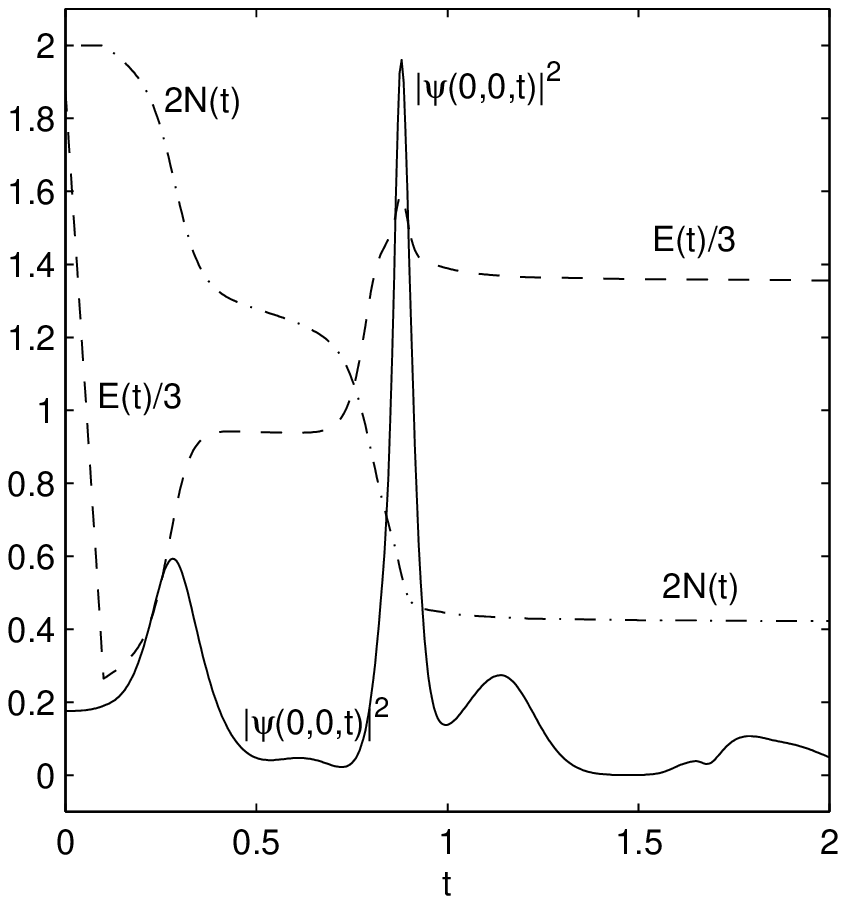,height=6cm,width=7cm,angle=0} 
\quad \psfig{figure=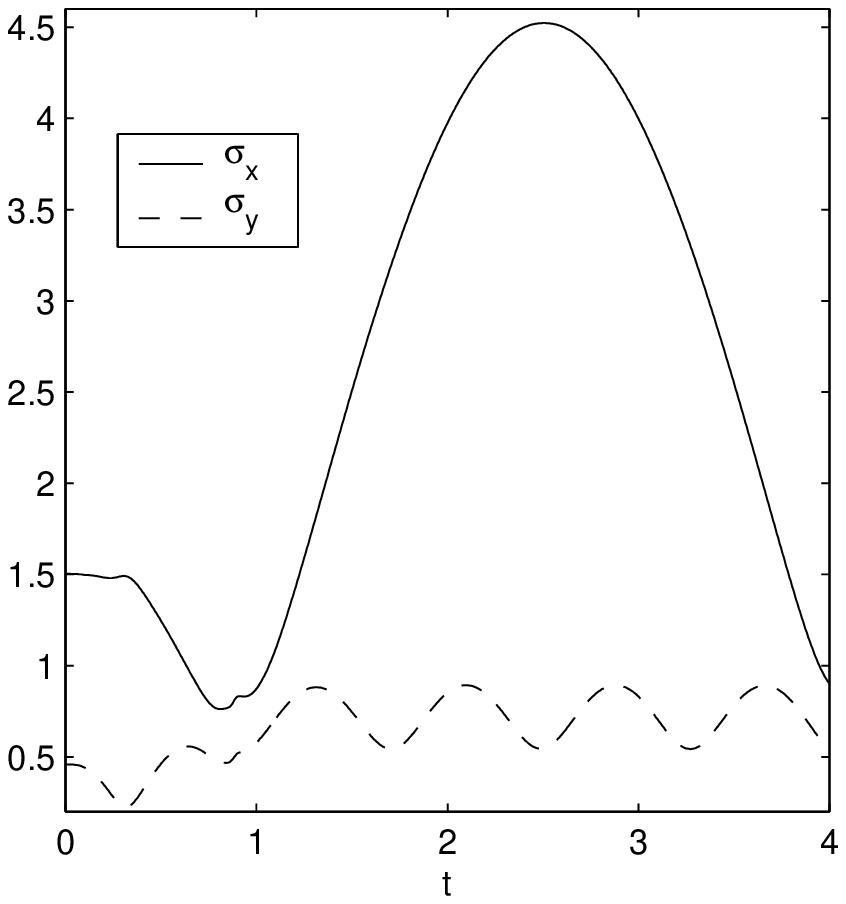,height=6cm,width=7cm,angle=0} }

Figure 7:  Numerical results  in Example 2 case III. \quad 
 a). Surface plot of the density $|\psi|^2$ 
 with $\dt=0.15$: At time $t=0.8$ (left column) and 
$t=3.2$ (right column). \quad 
 Normalization, energy and central density 
 $|\psi(0,0,t)|^2$ (left column) 
and condensate widths (right column) 
as functions of time:  b).  With $\dt=0.15$;  c). $\dt=0.005$. 
 
\end{figure} 
 
\begin{figure}[htb] 
\centerline{a).\psfig{figure=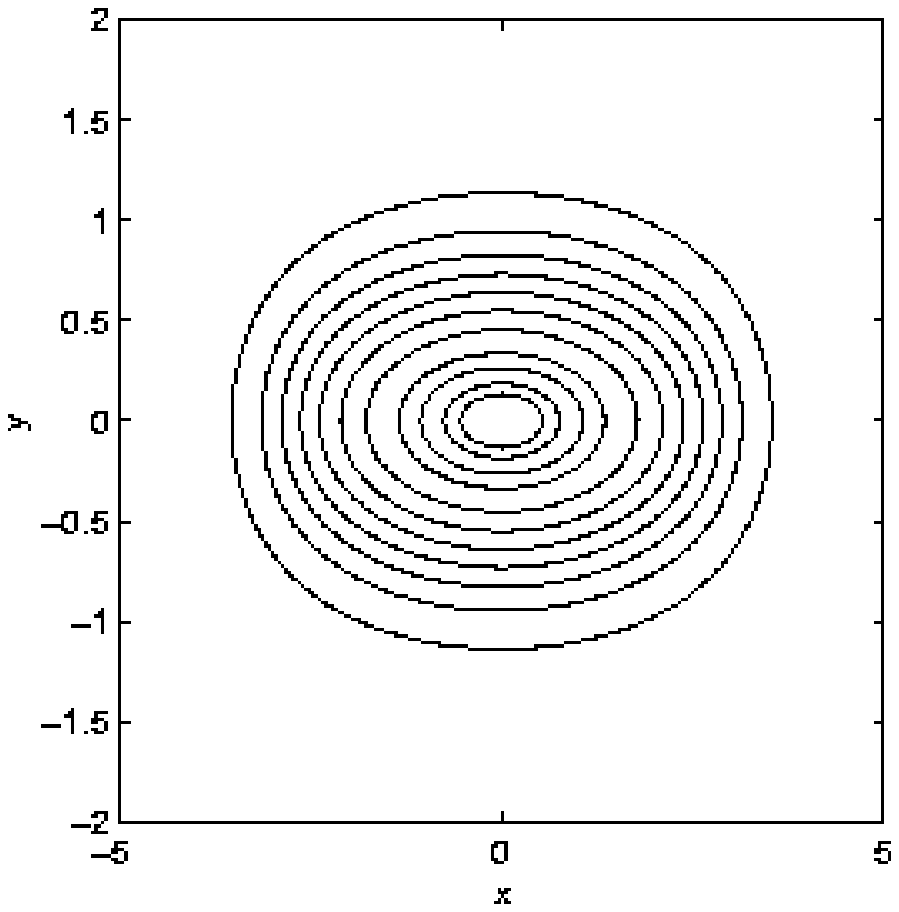,height=6cm,width=7cm,angle=0} 
\quad b)\psfig{figure=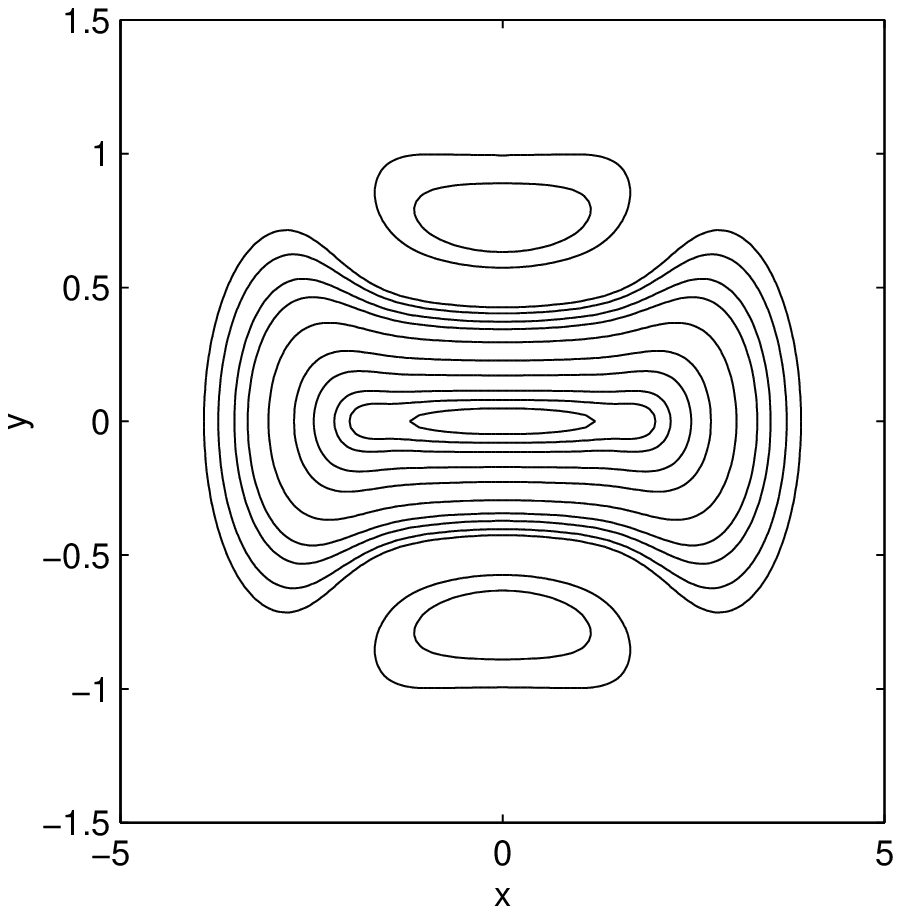,height=6cm,width=7cm,angle=0} } 
\centerline{c).\psfig{figure=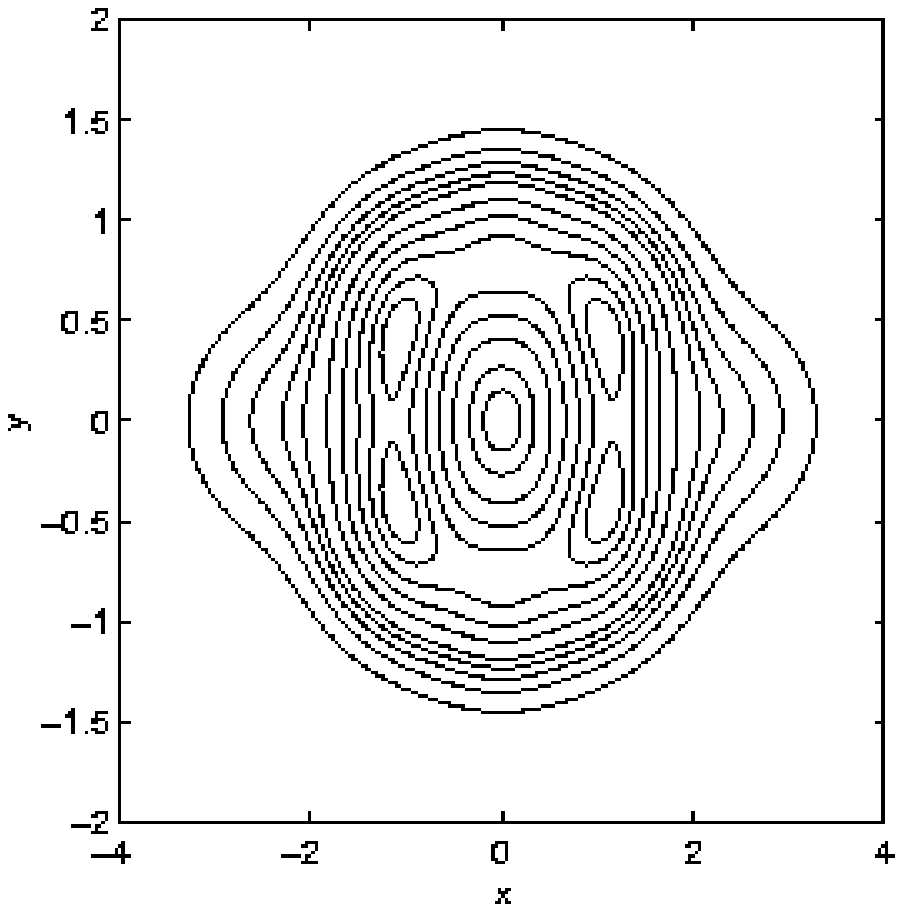,height=6cm,width=7cm,angle=0} 
\quad d)\psfig{figure=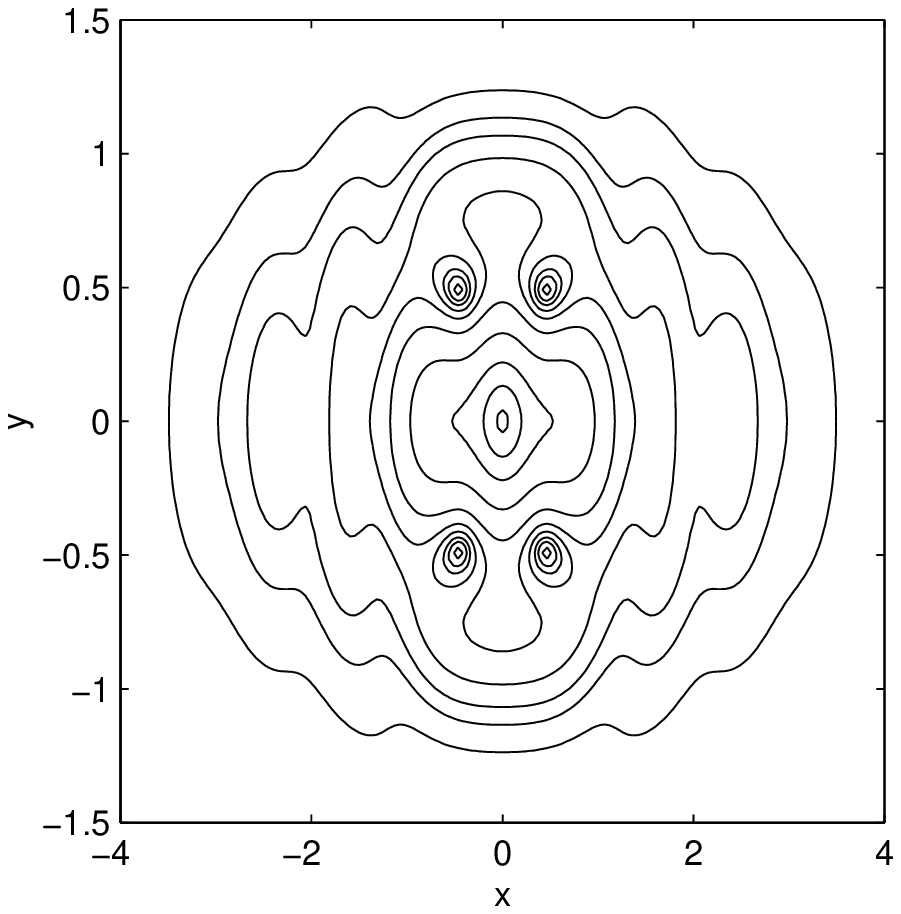,height=6cm,width=7cm,angle=0} } 
\centerline{e).\psfig{figure=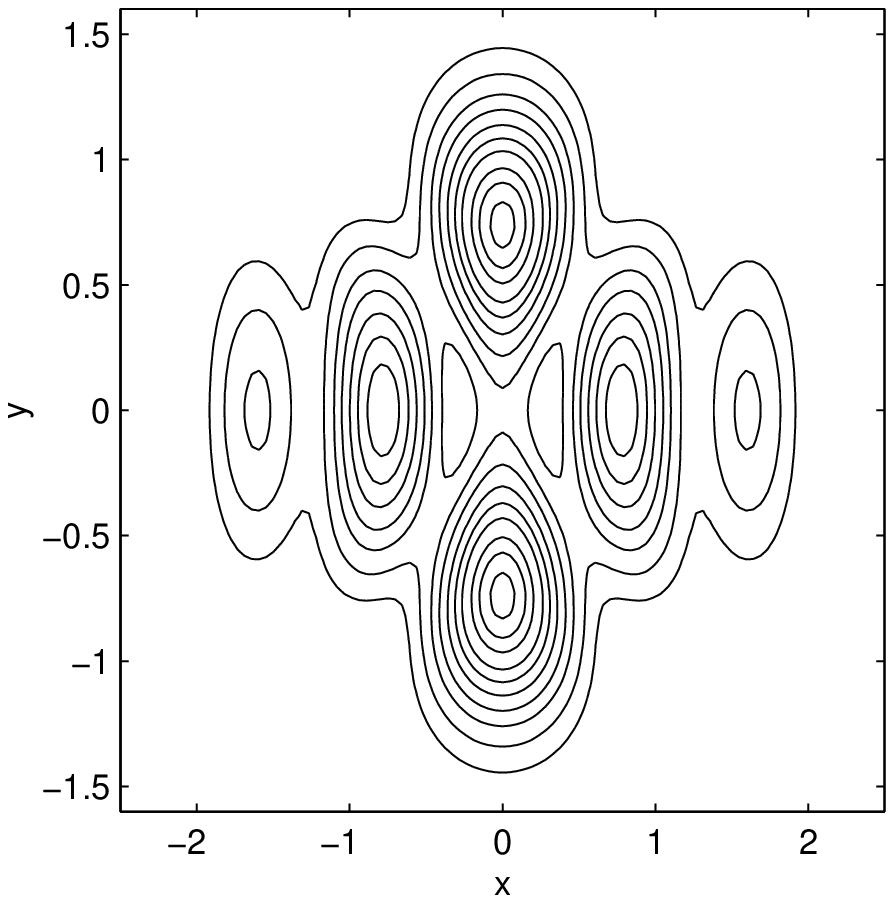,height=6cm,width=7cm,angle=0} 
\quad f)\psfig{figure=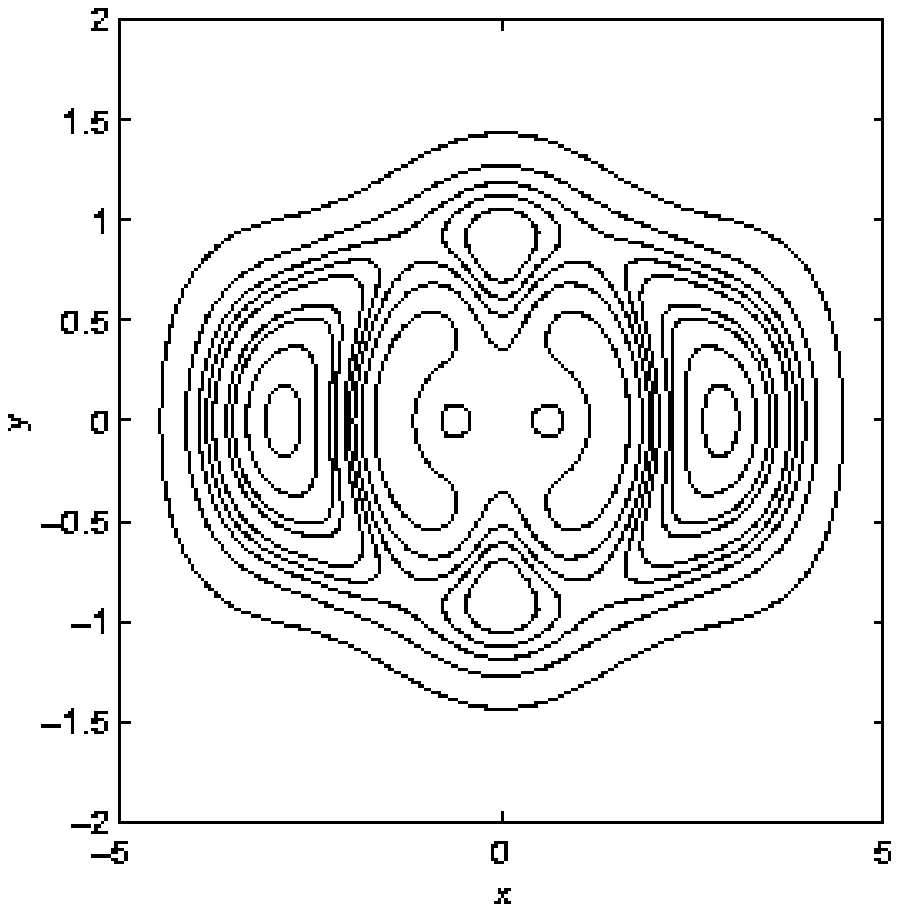,height=6cm,width=7cm,angle=0} } 
 
 Figure 8:  Contour plots of the density $|\psi|^2$ 
at different times in Example 2 case III with $\dt=0.15$. 
\  a). $t=0$, 
b). $t=0.4$, c). $t=0.8$, d). $t=1.2$, e). $t=1.6$, f). $t=2.4$. 
\end{figure} 
 
\end{document}